\documentclass[english]{elsarticle}
\usepackage[T1]{fontenc}
\usepackage[latin9]{inputenc}
\usepackage{geometry}
\usepackage{verbatim}
\usepackage{amsmath}
\usepackage{float}
\usepackage{amsbsy}
\usepackage{amssymb}
\usepackage{stackrel}
\usepackage{graphicx}
\usepackage{esint}
\usepackage{bm}
\usepackage{algpseudocode}
\usepackage{algorithm}
\usepackage{courier}

\usepackage{amsfonts} 

\usepackage{textcomp}
\usepackage{tabularx}
\usepackage{lineno}

\usepackage{todonotes}

\usepackage{enumitem}

\def\onedot{$\mathsurround0pt\ldotp$}

\def\cddot{
  \mathbin{\vcenter{\baselineskip .45ex
    \hbox{\onedot}\kern-.1ex
    \hbox{\onedot}}%
  }}%

\def\cdddot{
  \mathbin{\vcenter{\baselineskip .45ex
    \hbox{\onedot}\kern-.1ex
    \hbox{\onedot}\kern-.1ex
    \hbox{\onedot}%
  }}%
}

\usepackage[table,xcdraw]{xcolor}
\usepackage{adjustbox,lipsum}
\usepackage{subcaption}

\usepackage{xcolor}

\newcommand{\zero}{\bm{0}} 

\newcommand{\ap}[3]{{#1}_{#2}^{#3}} 

\newcommand{\vv}[3]{\mathbf{#1}_{#2}^{#3}} 
\newcommand{\vvm}[3]{\bm{#1}_{#2}^{#3}} 

\newcommand{\vvdot}[3]{\dot{\mathbf{#1}}_{#2}^{#3}} 
\newcommand{\vvmdot}[3]{\dot{\bm{#1}}_{#2}^{#3}} 

\newcommand{\vvddot}[3]{\Ddot{\mathbf{#1}}_{#2}^{#3}} 
\newcommand{\vvmddot}[3]{\Ddot{\bm{#1}}_{#2}^{#3}} 

\newcommand{\pder}[2]{\frac{\partial #1}{\partial #2}}

\newcommand{\etahat}{\hat{\eta}}

\newcommand{\R}[1]{\mathbb{R}^{#1}}

\newcommand{\dotpr}[1]{\cdot_{#1}}

\newcommand{\lman}{\vv{L}{\Gamma}{}}
\newcommand{\qman}{\vv{Q}{\Gamma}{}}
\newcommand{\lmane}{\vv{L}{\Gamma,e}{}}
\newcommand{\qmane}{\vv{Q}{\Gamma,e}{}}

\newcommand{\sscr}[1]{\textsuperscript{\textcolor{blue}{#1}}}

\newcommand{\newl}{\vspace{1\baselineskip} \par}



\usepackage{babel}

\journal{Computer Methods in Applied Mechanics and Engineering}

\begin{document}

\begin{frontmatter}{}

\title{Craig-Bampton-based Quadratic Manifold for Nonlinear Substructuring}

\author{Alexander Saccani\corref{cor1}}
\ead{asaccani@ethz.ch}

\author{Paolo Tiso}

\address{Institute for Mechanical Systems,\\ ETH Z\"urich, \\ Leonhardstrasse 21, 8092 Z\"urich, Switzerland}

\begin{abstract}
Component Mode Synthesis methods, such as the Craig-Bampton (CB) approach, are widely used in structural dynamics due to their modularity and compatibility with substructuring workflows. While highly effective for linear systems, extending these methods to geometrically nonlinear structures remains a significant challenge. 
In this work, we propose a nonlinear extension of the CB method tailored to such contexts. The approach is based on the construction of a quadratic reduction manifold, derived via perturbation analysis, in which high-frequency fixed-interface modes are statically condensed onto a reduced set of low-frequency modes and interface coordinates. This formulation enables the representation of geometric nonlinear effects without increasing the number of reduced degrees of freedom.
The resulting Nonlinear Craig-Bampton (NL-CB) reduced-order model is obtained through Galerkin projection onto the tangent space of the manifold and admits a polynomial structure that is efficient for time integration. The formulation preserves the Lagrangian structure of the underlying finite element model, ensuring consistent energetic behavior and numerical stability.
The proposed method is demonstrated on representative nonlinear structural systems of increasing complexity. The results show that the NL-CB model captures the essential nonlinear dynamic response while retaining the modularity and computational efficiency of classical substructuring approaches.

\end{abstract}

\begin{keyword}
Component Mode Synthesis \sep Dynamic Substructuring \sep Craig-Bampton \sep Geometric Nonlinearities \sep Quadratic Manifold \sep Reduced Order Model  
\end{keyword}

\end{frontmatter}{}

\section{Introduction} \label{sec:Introduction}
 The high computational cost associated with the dynamic analysis of structures using the Finite Element (FE) method led to the development of \textit{Component Mode Synthesis} techniques.
These procedures consist of partitioning large dimensional FE assemblies into smaller subdomains, referred to as substructures, computing a Reduced Order Model (ROM) for each substructure, and eventually assembling the substructure ROMs to retrieve a representative model of the global structural behavior.
By doing so, the computational effort for dynamic analysis is significantly reduced. 
\newl
The first pioneering work on CMS was independently developed by Guyan in \cite{Guyan1965} and Irons in \cite{Irons1965}.
The reduction procedure consisted of splitting the degrees of freedom (DoFs) of each substructure into few interface and more internal DOFs and in statically enslaving the internal DoFs to the interface DoFs. Thus, the substructure ROM was obtained by Galerkin projection of the FE substructure equations on a Reduction Basis (RB) of Static Modes (SMs).  This reduction technique goes by the name of the \textit{Guyan-Irons} method.
Hurty in \cite{Hurty1965} and Craig and Bampton in \cite{bampton1968} improved the accuracy of the Guyan-Irons model by complementing the RB with few low frequency \textit{Fixed Interface Modes} (FIMs).
This substructuring approach is known as the \textit{Craig-Bampton} method and has enjoyed successful application over the past decades.
As of today, the CB method shares its success in terms of popularity with the \textit{Rubin} substructuring method, first presented in \cite{Rubin1975}. The Rubin method can be regarded as the dual alternative to the CB method, in that \textit{Free-interface} Vibration Modes (VMs) and static solutions to imposed loads are used for reduction.
\newl
Substructuring techniques are well established for the reduction of linear systems and are routinely used in industry. However, their nonlinear counterpart is still in its infancy. 
Investigating the response of nonlinear structures plays a pivotal role in the design of a broad array of structures, such as aerospace components \cite{Deaton2015,Perez2014}, energy harvesters \cite{Fan2024}, MEMs \cite{Jia2012}, and friction joints  \cite{Segalman2006}, to name a few. 
Among different nonlinearity types, geometric nonlinearities are extremely relevant in the design of lightweight structures or systems featuring large displacements, as the in-service motion breaks the linearity assumption.
As a result, the FE equations are nonlinear. This implies a substantial increase of the computational cost compared to linear models.
Even more than for linear cases, nonlinear ROMs are usually the only feasible option for design optimization and parametric explorations.
\newl
Unlike dynamic substructuring approaches, monolithic model reduction techniques are well consolidated for geometrically nonlinear structures. 
In this field, popular approaches are \textit{Implicit Condensation} \cite{McEwan2001}, \textit{Implicit Condensation and Expansion} (ICE)\cite{Hollkamp2008}, the method of \textit{Static Modal Derivatives} (SMDs) \cite{Idelsohn1985,KaramoozMahdiabadi2021}, the method of \textit{Dual Modes} \cite{Kim2013}, the method of \textit{Quadratic Manifold} \cite{jain2017,Rutzmoser2017}, and the more recent method of \textit{Spectral Submanifolds} \cite{jain2022compute}\cite{Vizzaccaro2021}.
All of these methods can be viewed as an extension of modal analysis to nonlinear systems.
In fact, the ROM is constructed starting from a truncated set of VMs, either in the physical space or in state space.
The strength of these methods is that, in all cases, full order simulations are avoided, resulting in efficient ROM construction. 
\newl
In contrast to substructuring reduction, in monolithic approaches, the ROM is constructed starting from the assembled structure.
As a result, these techniques are not modular and therefore harder to integrate within industrial design procedures, where different structural components are designed independently by different teams.
Additionally, in case of repeated re-analysis of a structure with local substructure modifications, monolithic reduction forces a re-computation of the entire ROM.
Instead, if substructuring is used, when one component is modified, only the substructure ROM associated with that component is recomputed, while the other ROMs remain untouched. 
In this way, model updating is extremely efficient.
Furthermore, parallel computing can be exploited both to speed-up model construction and model evaluation.
Another valuable aspect of substructuring is the integration of equation-based models with experimental data, allowing for hybrid models that offer more flexibility in presence of modeling uncertainties \cite{Allen2019Substructuring}.
\newl
For the aforementioned reasons, an extension of CMS methods to geometrically nonlinear structures is extremely relevant for accelerating design in an industrial setting. 
To the best of the authors' knowledge, two main strategies have been probed so far towards the integration of monolithic ROM techniques with substructuring.
In the first strategy, the substructure ROMs are constructed in the same spirit of the ICE method, as first presented in \cite{Kuether2016}.
There, the CB RB was employed for substructure reduction and a third order polynomial restoring forces in both the FIMs and SMs amplitudes were identified from a set of nonlinear static solutions. 
This approach was integrated with \textit{Interface Reduction} methods \cite{Krattiger2019} in \cite{Kuether2017}, an essential step to obtain efficient substructuring ROMs.
Subsequently, the method was extended in \cite{KaramoozMahdiabadi2019} to the Rubin RB.
One of the main benefits of the ICE-substructuring approach lies in its non-intrusiveness, as the computation of the substructuring ROM does not require access to the FE source code.
However, as shown in \cite{Hippold2025}, the choice of the training load set strongly impacts the accuracy of the ensuing ROM. As a consequence, fine tuning of the imposed static loads is required.  
This limitation was also observed in ICE for monolithic reduction, as reported in \cite{Park2023,Shen2021}.
Compared to standard ICE, identifying the nonlinear restoring force for the substructuring ROM is even more challenging as loads imposed in the direction of SMs can trigger static buckling \cite{Hippold2025}. 
\newl
A second strategy for nonlinear substructuring consists of enriching the projection space used in CB and Rubin reduction with second order perturbation vectors, derived similarly to the SMDs.
This second approach was first used in \cite{Wenneker}, formalized in \cite{Wu2016} for the linear CB RB, and extended in \cite{Wu2019} to the Rubin RB.
Furthermore, more recent work in \cite{Mashayekhi2024, Bui2024} confirms the robustness of this approach.
Nevertheless, its efficiency is severely undermined by the quadratic growth of the number of complementary vectors with the size of the starting linear RB. 
For this reason, heuristic SMDs selection algorithms, as the one presented in \cite{Tiso2011}, have been used in \cite{Wu2016} to reduce the total number of companion vectors. Noteworthy, compared to the ICE substructuring approach, this last method is intrusive and requires the access to the FE source code to construct the ROM.
\newl
The work presented here can be regarded as modification of the CB substructuring ROM based on second order perturbation vectors in \cite{Wu2016}.
However, differently from \cite{Wu2016}, we use a quadratic manifold, in the same spirit of the one proposed in \cite{jain2017}, to avoid additional DOFs needed to capture the dynamic effects induced by the geometric nonlinearities.
Specifically, the reduction manifold is derived using perturbation theory,  enforcing a static condensation in modal coordinates of high frequency CB FIMs to low frequency CB FIMs and SMs.
The fundamental hypothesis behind the manifold derivation is that the inertia of high frequency FIMs can be neglected.
Starting from this hypothesis, we retrieve at linear level the classic CB ROM. Thus, the ROM presented here can be seen as an extension of the classic CB ROM to geometrically nonlinear structures, and, as such, we refer to it as NL-CB ROM.
The presented NL-CB ROM, derived upon Galerkin projection on the manifold tangent space, is in polynomial form, thus allowing remarkable speed-ups compared to the full FE model.
Additionally, we show that the NL-CB ROM preserves the \textit{Lagrangian-Structure} associated with the FE model \cite{Farhat2015}, thus guaranteeing good stability properties for time-integration.
\newl
The manuscript is divided into three main sections. After introducing the FE discretized model and briefly recalling the standard CB method in Section \ref{sec:preliminaries}, we present in Section \ref{sec:methodd} the NL-CB ROM. Eventually, in Section \ref{sec:applications}, we present numerical experiments where we apply the presented method on test cases of growing complexity.
\section{Set up} 

\label{sec:preliminaries}

\subsection{Full Order Model Equations}
\subsubsection{Substructures Dynamic Equilibrium}
We consider a generic partition of a FE assembly in $N_s$ different substructures.
For each of the substructures, $s \in \{{1,...,N_s}\}$, let $\vv{d}{}{(s)} = [\vv{u}{}{(s)T},\vv{q}{}{(s)T}]^T \in \R{n^{(s)}}$ be the vector of substructure DoFs, partitioned into internal DoFs, $\vv{u}{}{(s)} \in \R{n_u^{(s)}}$,  and interface DoFs, $\vv{q}{}{(s)} \in \R{n_q^{(s)}}$.
The dynamic equilibrium equations for a generic substructure $s$, assuming large displacements and rotations, write
\begin{equation}
    \vv{M}{}{(s)} \vvddot{d}{}{(s)} + \vv{D}{}{(s)} \vvdot{d}{}{(s)} + \vv{K}{}{(s)} \vv{d}{}{(s)} + \vv{f}{}{(s)}(\vv{d}{}{(s)}) = \vv{f}{}{ext(s)} + \vv{g}{}{(s)},
    \label{eq:subsEq}
\end{equation}
or, in partitioned format
\begin{equation}
\label{eq:subsEq1}
\begin{bmatrix}
    \vv{M}{uu}{(s)} & \vv{M}{uq}{(s)} \\
   \vv{M}{qu}{(s)} & \vv{M}{qq}{(s)} 
\end{bmatrix}
    \begin{bmatrix}
    \vvddot{u}{}{(s)} \\ \vvddot{q}{}{(s)} 
\end{bmatrix}
+
\begin{bmatrix}
    \vv{D}{uu}{(s)} & \vv{D}{uq}{(s)} \\
   \vv{D}{qu}{(s)} & \vv{D}{qq}{(s)} 
\end{bmatrix}
    \begin{bmatrix}
    \vvdot{u}{}{(s)} \\ \vvdot{q}{}{(s)} 
\end{bmatrix}
+
\begin{bmatrix}
    \vv{K}{uu}{(s)} & \vv{K}{uq}{(s)} \\
   \vv{K}{qu}{(s)} & \vv{K}{qq}{(s)} 
\end{bmatrix}
\begin{bmatrix}
    \vv{u}{}{(s)} \\ \vv{q}{}{(s)} 
\end{bmatrix}
+
\begin{bmatrix}
    \vv{f}{u}{(s)}(\vv{u}{}{(s)},\vv{q}{}{(s)}) \\ 
    \vv{f}{q}{(s)}(\vv{u}{}{(s)},\vv{q}{}{(s)})
\end{bmatrix}
=
\begin{bmatrix}
    \vv{f}{u}{ext(s)}(t) \\ \vv{f}{q}{ext(s)}(t)
\end{bmatrix}
+
\begin{bmatrix}
    \zero \\ \vv{g}{q}{(s)}
\end{bmatrix}.
\end{equation}
Here, $\vv{M}{}{(s)}$, $\vv{K}{}{(s)}$, and $\vv{D}{}{(s)}$ $\in \R{n^{(s)} \times n^{(s)}}$ are the mass, damping, and stiffness matrix; $\vv{f}{}{(s)} \in \R{n^{(s)}}$ is the vector of nonlinear elastic forces; $\vv{f}{}{ext(s)} \in \R{n^{(s)}}$ the vector of external time-dependent loads, and $\vv{g}{}{(s)} \in \R{n^{(s)}}$ the vector of reaction interface forces exerted by the connected substructures.\\
Furthermore, when linear elastic material constitutive law is used, the vector of elastic forces originating from FE formulations based on the \textit{Green-Lagrange} or \textit{von K\'arm\'an} strain measures \cite{crisfield1991} are third order polynomials in the nodal variables \cite{Mignolet2008}, writing
\begin{equation}
\label{eq:elastForcPolSubs}
    \vv{f}{}{(s)}(\vv{d}{}{(s)}) = \vv{K}{}{2(s)} : (\vv{d}{}{(s)} \otimes \vv{d}{}{(s)}) + \vv{K}{}{3(s)} \cdddot (\vv{d}{}{(s)} \otimes \vv{d}{}{(s)} \otimes \vv{d}{}{(s)}),
\end{equation}
or, in partitioned format,
\begin{subequations}
\label{eq:NlTensSubs}
    \begin{align}
   \vv{f}{u}{}(\vv{u}{}{},\vv{q}{}{}) =&\ \vv{K}{uuu}{2}:\vv{u}{}{}\otimes \vv{u}{}{} + \vv{K}{uuq}{2}:\vv{u}{}{}\otimes \vv{q}{}{} + \vv{K}{uqq}{2}:\vv{q}{}{}\otimes \vv{q}{}{} + \vv{K}{uuuu}{3} \cdddot \ \vv{u}{}{}\otimes \vv{u}{}{} \otimes \vv{u}{}{} + \nonumber \\
        &+ \vv{K}{uuuq}{3} \cdddot \ \vv{u}{}{}\otimes \vv{u}{}{} \otimes \vv{q}{}{} + \vv{K}{uuqq}{3} \cdddot \ \vv{u}{}{}\otimes \vv{q}{}{} \otimes \vv{q}{}{} +\vv{K}{uqqq}{3} \cdddot \ \vv{q}{}{}\otimes \vv{q}{}{} \otimes \vv{q}{}{},\\
\vv{f}{q}{}(\vv{u}{}{},\vv{q}{}{}) =&\ \vv{K}{quu}{2}:\vv{u}{}{}\otimes \vv{u}{}{} + \vv{K}{quq}{2}:\vv{u}{}{}\otimes \vv{q}{}{} + \vv{K}{qqq}{2}:\vv{q}{}{}\otimes \vv{q}{}{} + \vv{K}{quuu}{3} \cdddot \ \vv{u}{}{}\otimes \vv{u}{}{} \otimes \vv{u}{}{} + \nonumber \\
        &+ \vv{K}{quuq}{3} \cdddot \ \vv{u}{}{}\otimes \vv{u}{}{} \otimes \vv{q}{}{} + \vv{K}{quqq}{3} \cdddot \ \vv{u}{}{}\otimes \vv{q}{}{} \otimes \vv{q}{}{} +\vv{K}{qqqq}{3} \cdddot \ \vv{q}{}{}\otimes \vv{q}{}{} \otimes \vv{q}{}{}.
\end{align}
\end{subequations}
Here, we omit the superscript "$(s)$" to maintain the notation compact.
In this expression, "$\otimes$","$:$", and "$\cdddot$" denote, respectively, the dyadic product, double, and triple tensor contraction; $\vv{K}{}{2}$ and $\vv{K}{}{3}$ are  the tensors associated with the quadratic and cubic components of elastic forces.
Since the elastic forces are derived through differentiation of an elastic potential energy, the linear stiffness matrix, $\vv{K}{}{(s)}$, and the nonlinear stiffness tensors, $\vv{K}{}{2(s)}$ and $\vv{K}{}{3(s)}$, are symmetric with respect to all axes \cite{rutzmoserThesis}. 
Additionally, these operators are extremely sparse, which is an essential feature for efficient storage \cite{rutzmoserThesis}.
\subsubsection{Assembled FE Equations}
\label{sec:assembledFEM}
Compatibility conditions on interface displacements and forces must be imposed to solve the dynamic problem.
Displacement compatibility is written as
\begin{equation}
   \sum_{i=1}^{N_s} \vv{B}{}{(s)}\vv{d}{}{(s)} = \zero.
   \label{eq:dispCompEq}
\end{equation}
where $\vv{B}{}{(s)} \in \R{p\times n^{(s)}}$.
This condition enforces continuity of the displacement field at the interface between two or more substructures.
Conversely, equilibrium at interfaces writes
\begin{equation}
    \vv{g}{}{(s)} = \vv{B}{}{(s)T}\vvm{\lambda}{}{},
    \label{eq:forcCompEq}
\end{equation}
where  $\vvm{\lambda}{}{} \in \R{p}$, the generalized interface force vector, is required to prevent tearing of the substructures at the interface.\newl
The equations for the assembled structures are obtained following a \textit{Primal Assembly} procedure, where the interface reaction forces are eliminated by projecting the substructure equations on a basis of compatible displacements.
This operation is similar to the standard assembly procedure in the FE method, where the substructures play the role of elements.
Specifically, we first define a vector $\vv{\tilde{d}}{}{}\in \R{\tilde{n}}$ that collects all the displacement variables from the different substructures, including only once interface displacements shared between different substructures.
In this way, the displacement variables of substructure $s$ can be extracted as
\begin{equation}
    \vv{d}{}{(s)} = \vv{L}{}{(s)} \vv{\tilde{d}}{}{},
    \label{eq:subsLmatr}
\end{equation}
where $\vv{L}{}{(s)} \in \R{n^{(s)\times\tilde{n}}}$ is a boolean extraction matrix.
By inserting Eq. \eqref{eq:subsLmatr} in the substructure dynamic equilibrium in Eq. \eqref{eq:subsEq}, projecting on $\vv{L}{}{(s)}$, and summing over $s$, we retrieve the assembled dynamic equilibrium
\begin{equation}
    \vv{M}{}{} \Ddot{\tilde{\mathbf{d}}} + \vv{D}{}{} \dot{\tilde{\mathbf{d}}} + \vv{K}{}{} \vv{\tilde{d}}{}{} +\vv{f}{}{}(\vv{\tilde{d}}{}{}) = \vv{f}{}{ext}(t),
\end{equation}
where
\begin{subequations}
\begin{gather}
     \vv{M}{}{} = \sum_{i=1}^{N_s} \vv{L}{}{(s)T} \vv{M}{}{(s)} \vv{L}{}{(s)},\quad  \vv{D}{}{} =  \sum_{i=1}^{N_s} \vv{L}{}{(s)T} \vv{D}{}{(s)} \vv{L}{}{(s)}, \quad \vv{K}{}{} = \sum_{i=1}^{N_s}  \vv{L}{}{(s)T} \vv{K}{}{(s)} \vv{L}{}{(s)},\\  \quad  \vv{f}{}{}(\vv{\tilde{d}}{}{}) = \sum_{i=1}^{N_s}  \vv{L}{}{(s)T}\vv{f}{}{(s)}( \vv{L}{}{(s)}\vv{\tilde{d}}{}{(s)}),\quad \vv{f}{}{ext}(t) = \sum_{i=1}^{N_s} \vv{L}{}{(s)T} \vv{f}{}{ext(s)}(t).
\end{gather}
\end{subequations}
In the previous expression, the projected interface forces vanish in the summation, since $\vv{L}{}{(s)}$ satisfy
\begin{equation}
    \sum_{i=1}^{N_s} \vv{B}{}{(s)}\vv{L}{}{(s)} = \zero.
    \label{eq:condCompDispB}
\end{equation}
\subsection{Craig-Bampton Reduction for Linear Systems}
The CB model reduction technique, originally presented in the seminal paper in \cite{craig1968coupling} to reduce linear structural systems, consists of two subsequent steps:
a reduction at substructure level using the Galerkin projection method and a subsequent primal assembly of the reduced substructure ROMs.
Specifically, the DoFs of substructure $s$ are approximated as
\begin{equation}
    \vv{d}{}{(s)} \approx \vv{V}{}{(s)}\vvm{\xi}{}{(s)},
    \label{eq:CBbasis1}
\end{equation}
where $\vvm{\xi}{}{(s)} \in \R{ \ap{m}{}{(s)}}$, $\ap{m}{}{(s)} \ll \ap{n}{}{(s)}$, is the vector of substructure reduced coordinates and $\vv{V}{}{(s)} \in \R{\ap{n}{}{(s)} \times \ap{m}{}{(s)}}$ is the substructure RB.
Splitting the substructure  displacement vector into internal and interface DoFs, the CB approximation in Eq. \eqref{eq:CBbasis1} writes
\begin{equation}
\begin{bmatrix}
        \vv{u}{}{(s)}\\
        \vv{q}{}{(s)}
    \end{bmatrix} \approx
    \begin{bmatrix}
        \vvm{\Phi}{}{(s)} & \vv{S}{}{(s)} \\ \zero & \mathbf{I}
    \end{bmatrix}
    \begin{bmatrix}
        \vvm{\eta}{}{(s)}\\
        \vv{q}{}{(s)}
    \end{bmatrix}
    \label{eq:CBbasis2}
\end{equation}
Here, $\vvm{\Phi}{}{(s)} \in \R{\ap{n}{u}{(s)} \times \ap{n}{\phi}{(s)}}$ contains few low frequency FIMs, while $\vv{S}{}{(s)} \in \R{\ap{n}{u}{(s)} \times \ap{n}{q}{(s)}}$ gathers the so called SMs, equivalently known as \textit{Component Modes}.
The FIMs basis is defined as
\begin{equation}
    \vvm{\Phi}{}{(s)} = [...,\vvm{\phi}{j}{(s)},...], \quad \text{with} \quad (-\ap{\omega}{j}{(s)2}\vv{M}{uu}{(s)} + \vv{K}{uu}{(s)})\vvm{\phi}{j}{(s)} = \zero,
\label{eq:fixedInterfModes}
\end{equation}
where $\vvm{\phi}{j}{(s)} \in \R{\ap{n}{u}{(s)}}$ is the $j$th VM of the substructure when all the substructure interface DoFs are constrained to zero, and $\ap{\omega}{j}{(s)}$ its associated angular frequency.
Additionally, the SMs basis, defined as
\begin{equation}
    \vv{S}{}{(s)} = - \vv{K}{uu}{(s)-1}\vv{K}{uq}{(s)},
    \label{eq:stModesBasis}
\end{equation}
represents the set of static displacements in the internal DoFs resulting from imposing unitary interface displacements.
The SMs are essential to capture the contributions of high frequency FIMs to the response to the elastic forces arising from interface motion.
\newl
In some applications, a further reduction of the substructure DoFs is achieved by resorting to interface reduction techniques \cite{Krattiger2019, Castanier2001, Hong2013, Lindberg2013}.
Specifically, the interface DoFs are reduced with the approximation 
\begin{equation}
    \vv{q}{}{(s)} \approx \vvm{\Psi}{}{(s)}\vvm{\chi}{}{(s)},
    \label{eq:intRedBasis}
\end{equation}
where $\vvm{\Psi}{}{(s)} \in \R{\ap{n}{q}{(s)} \times \ap{n}{\chi}{(s)} }$, with $\ap{n}{\chi}{(s)} \ll \ap{n}{q}{(s)}$, is a RB of \textit{Interface Modes}  and $\vvm{\chi}{}{(s)} \in \R{\ap{n}{\chi}{(s)}}$ is the vector of interface reduced coordinates.
Thus, the CB approximation in Eq. \eqref{eq:CBbasis2} combined with interface reduction becomes
\begin{equation}
\begin{bmatrix}
        \vv{u}{}{(s)}\\
        \vv{q}{}{(s)}
    \end{bmatrix} \approx
    \begin{bmatrix}
        \vvm{\Phi}{}{(s)} & \vv{S}{}{(s)} \vvm{\Psi}{}{(s)} \\ \zero & \vvm{\Psi}{}{(s)}
    \end{bmatrix}
    \begin{bmatrix}
        \vvm{\eta}{}{(s)}\\
        \vvm{\chi}{}{(s)}.
    \end{bmatrix}
    \label{eq:CBbasisIntRed}
\end{equation}
\newl
After defining the RB, the substructure reduced equations are obtained by inserting the approximation in Eq. \eqref{eq:CBbasis1} into the substructure dynamic equilibrium in Eq. \eqref{eq:subsEq} - neglecting nonlinear forces -  and by performing a Galerkin projection, thus obtaining
\begin{equation}
    \vv{M}{r}{(s)}\vvmddot{\xi}{}{(s)} + \vv{D}{r}{(s)} \vvmdot{\xi}{}{(s)} +\vv{K}{r}{(s)} \vvm{\xi}{}{(s)} = \vv{f}{r}{ext(s)}(t) + \vv{V}{}{(s)T}\vv{g}{}{(s)},
    \label{eq:CBsubsProj}
\end{equation}
where
\begin{equation}
    \vv{M}{r}{(s)} = \vv{V}{}{(s)T} \vv{M}{}{(s)} \vv{V}{}{(s)}, \quad 
    \vv{D}{r}{(s)} = \vv{V}{}{(s)T} \vv{D}{}{(s)} \vv{V}{}{(s)}, \quad
    \vv{K}{r}{(s)} = \vv{V}{}{(s)T} \vv{K}{}{(s)} \vv{V}{}{(s)}, \quad
    \vv{f}{r}{ext(s)} = \vv{V}{}{(s)T}  \vv{f}{}{ext(s)}(t).
\end{equation}
Next, similarly to what described in Section \ref{sec:assembledFEM}, the substructure ROMs are assembled using a primal assembly approach.
After defining the global vector of reduced displacements, $\vvm{\tilde{\xi}}{}{}\in\R{\tilde{m}}$, and the localization matrices  $\vv{L}{r}{(s)} \in \R{m^{(s)}\times \tilde{m}}$, such that 
\begin{equation}
    \vvm{\xi}{}{(s)} = \vv{L}{r}{(s)} \vvm{\tilde{\xi}}{}{},
    \label{eq:extrRed}
\end{equation}
the assembled CB ROM equations writes\
\begin{equation}
        \vv{M}{r}{}\vvmddot{\tilde{\xi}}{}{} + \vv{D}{r}{} \vvmdot{\tilde{\xi}}{}{} +\vv{K}{r}{} \vvm{\tilde{\xi}}{}{} = \vv{f}{r}{ext}(t),
        \label{eq:CBfinal}
\end{equation}
where
\begin{equation}
    \vv{M}{r}{} = \sum_{s=1}^{N_s} \vv{L}{r}{(s)T} \vv{M}{r}{(s)} \vv{L}{r}{(s)}, \quad 
    \vv{D}{r}{} = \sum_{s=1}^{N_s} \vv{L}{r}{(s)T} \vv{D}{r}{(s)} \vv{L}{r}{(s)}, \quad
    \vv{K}{r}{} = \sum_{s=1}^{N_s} \vv{L}{r}{(s)T} \vv{K}{r}{(s)} \vv{L}{r}{(s)}, \quad
    \vv{f}{r}{ext} = \sum_{s=1}^{N_s} \vv{L}{r}{(s)T}  \vv{f}{r}{ext(s)}(t).
\end{equation}

\section{Nonlinear Craig-Bampton Reduction}
\label{sec:methodd}
In this Section we present our generalization of CB reduction to geometrically nonlinear structures, referred to as NL-CB ROM.
The NL-CB ROM is constructed by reducing the substructure dynamic equilibrium equations using a Galerkin projection on a nonlinear second order manifold.
We illustrate the theory in three different subsections: in Subsection \ref{sec:manDer} we derive the theoretical basis underlying the substructure reduction manifold;  in Subsection \ref{sec:manComp} we show how the substructure manifold and the substructure reduced equations are numerically computed; finally, in Subsection \ref{sec:romass} we detail the assembly of the NL-CB ROM using the primal approach.
\subsection{Substructure Reduction}
\label{sec:method}
\subsubsection{Substructure Manifold Derivation}
\label{sec:manDer}
 The NL-CB ROM reduction subspace is derived assuming an unforced structure.
This choice does not restrict the application of the resulting ROM only to the investigation of the unforced response, allowing for forced predictions. 
In fact, the use of reduction subspaces derived assuming unforced dynamics for forced  predictions is at the heart of numerous successful model reduction methods, such as Modal Analysis \cite{geradin2015mechanical}, ICE \cite{Hollkamp2008}, the methods of Modal Derivatives \cite{Idelsohn1985}, Nonlinear Normal Modes \cite{Kerschen2009}, and Spectral Submanifolds \cite{haller2016,jain2022compute,Vizzaccaro2021}.\newl
We consider the unforced dynamic equilibrium of the internal DoFs of a generic substructure $s$  -- see first block row of  Eq. \eqref{eq:subsEq1}.
Furthermore, we neglect damping forces as well as the coupling inertial terms between interface and internal DoFs, as these terms are usually small compared to cross coupling elastic forces.
Thus, the dynamic equilibrium writes
\begin{equation}
    \vv{M}{uu}{} \vvddot{u}{}{} + \vv{K}{uu}{} \vv{u}{}{} + \vv{K}{uq}{} \vv{q}{}{} + \vv{f}{u}{}(\vv{u}{}{},\vv{q}{}{}) = \zero.
    \label{eq:eqIntApp}
\end{equation}
Here, as in the rest of this section, we omit the superscript $"(s)"$ when referring to substructure related quantities and reinstate it only when strictly necessary.
Although derived for zero external forcing, Eq.\eqref{eq:eqIntApp} describes the  motion of internal DoFs $\vv{u}{}{}$, forced by the elastic forces $\vv{f}{u}{}$ resulting from the motion of interface DoFs $\vv{q}{}{}$.\newl
To obtain a minimal dimensional ROM, we derive the manifold equations allowing for interface reduction.
Thus, we use the approximation in Eq. \eqref{eq:intRedBasis} for interface DoFs, from which, however, the unreduced interface case can be retrieved using $\vvm{\Psi}{}{}$ equal to the identity matrix.
In this way, Eq. \eqref{eq:eqIntApp} becomes
\begin{equation}
       \vv{M}{uu}{} \vvddot{u}{}{} + \vv{K}{uu}{} \vv{u}{}{} + \vv{K}{uq}{} \vvm{\Psi}{}{} \vvm{\chi}{}{} + \vv{f}{u}{}(\vv{u}{}{},\vvm{\Psi}{}{}\vvm{\chi}{}{}) = \zero.
    \label{eq:eqIntApp1}
\end{equation}
Next, we perform a change of variables, writing the internal DoFs as a linear combination of FIMs  (defined in Eq.\eqref{eq:fixedInterfModes}), divided between low and high frequency modes.
Specifically, we write
\begin{equation}
    \vv{u}{}{} = [\vvm{\Phi}{}{},\vvm{\hat{\Phi}}{}{}]
    \begin{bmatrix}
        \vvm{\eta}{}{} \\ \vvm{\hat{\eta}}{}{}
    \end{bmatrix},
    \label{eq:changVar}
\end{equation}
where $\vvm{\Phi}{}{} \in \R{\ap{n}{u}{} \times \ap{n}{\phi}{}}$ collects the first $\ap{n}{\phi}{}$ FIMs, while $\vvm{\hat{\Phi}}{}{} \in \R{\ap{n}{u}{} \times (\ap{n}{u}{} - \ap{n}{\phi}{}) }$ the remaining $\ap{n}{\hat{\phi}}{} = \ap{n}{u}{} - \ap{n}{\phi}{}$ high frequency FIMs.
In the new coordinates, Eq. \eqref{eq:eqIntApp1} writes
\begin{align}
\label{eq:dynEqModal}
    \begin{bmatrix}
        \vvm{\Phi}{}{T}\vv{M}{uu}{}\vvm{\Phi}{}{} & \zero \\
        \zero & \vvm{\hat{\Phi}}{}{T}\vv{M}{uu}{}\vvm{\hat{\Phi}}{}{}
    \end{bmatrix}
    \begin{bmatrix}
        \vvmddot{\eta}{}{} \\ \vvmddot{\hat{\eta}}{}{}
    \end{bmatrix} +
    &
        \begin{bmatrix}
        \vvm{\Phi}{}{T}\vv{K}{uu}{}\vvm{\Phi}{}{} & \zero \\
        \zero & \vvm{\hat{\Phi}}{}{T}\vv{K}{uu}{}\vvm{\hat{\Phi}}{}{}
    \end{bmatrix}
    \begin{bmatrix}
        \vvm{\eta}{}{} \\ \vvm{\hat{\eta}}{}{}
    \end{bmatrix}+
    \begin{bmatrix}
        \vvm{\Phi}{}{T} \vv{K}{uq}{}\vvm{\Psi}{}{}\vvm{\chi}{}{} \\ 
        \vvm{\hat{\Phi}}{}{T} \vv{K}{uq}{}\vvm{\Psi}{}{}\vvm{\chi}{}{}
    \end{bmatrix} + \\ + &
    \begin{bmatrix}
        \vvm{\Phi}{}{T} \vv{f}{u}{} (  \vvm{\Phi}{}{}\vvm{\eta}{}{} + \vvm{\hat{\Phi}}{}{}\vvm{\hat{\eta}}{}{}, \vvm{\Psi}{}{}\vvm{\chi}{}{}) \\
        \vvm{\hat{\Phi}}{}{T} \vv{f}{u}{} (  \vvm{\Phi}{}{}\vvm{\eta}{}{} + \vvm{\hat{\Phi}}{}{}\vvm{\hat{\eta}}{}{}, \vvm{\Psi}{}{}\vvm{\chi}{}{})
    \end{bmatrix}
    =
    \begin{bmatrix}
        \zero \\
        \zero
    \end{bmatrix},
    \nonumber
\end{align}
where we use mass and stiffness orthogonality of FIMs \cite{geradin2015mechanical}.
Notably, the coupling in the equations stems from two terms: a linear term, which is a function of the interface DoFs only, and a nonlinear term, which depends on both interface DoFs and FIMs coordinates.
\newl
In this work, the reduction manifold for the substructure equation of motion is obtained assuming that high frequency FIMs are statically enslaved to low frequency FIMs and SMs.
Specifically, by neglecting the inertia forces in the second equation in \eqref{eq:dynEqModal}  we obtain
\begin{equation}
    \label{eq:staticEqModes}
    \vvm{\hat{\Phi}}{}{T}\vv{K}{uu}{}\vvm{\hat{\Phi}}{}{}\vvm{\hat{\eta}}{}{} + \vvm{\hat{\Phi}}{}{T} \vv{K}{uq}{}\vvm{\Psi}{}{}\vvm{\chi}{}{} +  \vvm{\hat{\Phi}}{}{T} \vv{f}{u}{} (  \vvm{\Phi}{}{}\vvm{\eta}{}{} + \vvm{\hat{\Phi}}{}{}\vvm{\hat{\eta}}{}{}, \vvm{\Psi}{}{}\vvm{\chi}{}{}) = \zero.
\end{equation}
We expect that this simplifying assumption holds, as we assume that the excitation to high frequency FIMs, coming from low frequency FIMs and interface motion, is at low frequency compared to their characteristic frequencies.
A similar assumption lies behind the monolithic reduction strategies underlying ICE \cite{Shen2021}, SMDs \cite{jain2017,KaramoozMahdiabadi2021,Weeger2016}, and Static Quadratic Manifold ROMs \cite{jain2017}.
\newl
Equation \eqref{eq:staticEqModes} implicitly defines a static manifold:
the amplitudes of high frequency FIMs can be expressed as a function of the amplitudes of low frequency FIMs and of SMs.
Formally, if we write Eq. \eqref{eq:staticEqModes} as
\begin{equation}
    \vvm{g}{}{}(\vvm{\hat{\eta}}{}{},\vvm{\eta}{}{},\vvm{\chi}{}{}) = \zero,
\end{equation}
 we can define the static manifold $\vvm{\Gamma}{\etahat}{}$,  in a neighborhood of the equilibrium point, as
\begin{equation}
   \vvm{\Gamma}{\etahat}{}:= \vvm{\eta}{}{} \times \vvm{\chi}{}{} \mapsto \vvm{\hat{\eta}}{}{} \quad : \quad \vvm{g}{}{}(\vvm{\Gamma}{\etahat}{}(\vvm{\eta}{}{},\vvm{\chi}{}{}),\vvm{\eta}{}{},\vvm{\chi}{}{}) = \zero.
   \label{eq:stMan}
\end{equation}
This statement can be justified using the Implicit Function Theorem \cite{Zeidler1986}, as the Jacobian of $\vvm{g}{}{}$ with respect to $\vvm{\hat{\eta}}{}{}$, evaluated at the origin is $\vvm{\hat{\Phi}}{}{T}\vv{K}{uu}{}\vvm{\hat{\Phi}}{}{}$, which is non singular.
Additionally, since $\vvm{g}{}{}$ is analytic in its arguments - as we work with polynomial nonlinearities in Eq. \eqref{eq:NlTensSubs} - we can expand the manifold in Eq. \eqref{eq:stMan} in a convergent Taylor series of $\vvm{\eta}{}{}$ and $\vvm{\chi}{}{}$, truncated at second order.
 In this way, we can write the manifold approximation as
\begin{equation}
    \vvm{\Gamma}{\etahat}{}(\vvm{\eta}{}{},\vvm{\chi}{}{}) \approx  \vv{\hat{B}}{}{}\vvm{\chi}{}{} +\vv{\hat{C}}{}{}:\vvm{\eta}{}{} \otimes \vvm{\eta}{}{}  +\vv{\hat{D}}{}{} : \vvm{\eta}{}{} \otimes \vvm{\chi}{}{} +\vv{\hat{E}}{}{}: \vvm{\chi}{}{} \otimes \vvm{\chi}{}{},
    \label{eq:manFast}
\end{equation}
where 
$\vv{\hat{B}}{}{} \in \R{\ap{n}{\hat{\phi}}{} \times \ap{n}{\chi}{}}$,
$\vv{\hat{C}}{}{} \in \R{\ap{n}{\hat{\phi}}{} \times \ap{n}{\phi}{}  \times \ap{n}{\phi}{} }$, $\vv{\hat{D}}{}{} \in \R{\ap{n}{\hat{\phi}}{} \times \ap{n}{\phi}{}  \times \ap{n}{\chi}{} }$, and
$\vv{\hat{E}}{}{} \in \R{\ap{n}{\hat{\phi}}{} \times \ap{n}{\chi}{}  \times \ap{n}{\chi}{} }$.
The analytical expressions for $\vv{\hat{B}}{}{},\vv{\hat{C}}{}{},\vv{\hat{D}}{}{}, \text{and}\  \vv{\hat{E}}{}{}$ are presented in \ref{app:secOrdS}. 
The derivation of the manifold coefficients is performed by leveraging the tensorial decomposition of elastic forces in Eq. \eqref{eq:NlTensSubs}.
Interestingly, the Taylor expansion does not contain any linear term in $\vvm{\eta}{}{}$: these are identically zero as a consequence of linear FIMs decoupling -- see \ref{app:secOrdS}. \newl
When we insert Eq. \eqref{eq:manFast} into Eq. \eqref{eq:changVar}, we obtain
\begin{subequations}
    \label{eq:manApprox}
 \begin{align}
    \label{eq:manApproxA}
    \vv{u}{}{} &\approx \vvm{\Phi}{}{}\vvm{\eta}{}{} + \vvm{\hat{\Phi}}{}{}\vvm{\Gamma}{\etahat}{}(\vvm{\eta}{}{},\vvm{\chi}{}{}) \\
    &  \approx \vvm{\Phi}{}{} \vvm{\eta}{}{} + \vvm{\hat{\Phi}}{}{} \vv{\hat{B}}{}{}\vvm{\chi}{}{} +  \vvm{\hat{\Phi}}{}{} \cdot \vv{\hat{C}}{}{}:\vvm{\eta}{}{} \otimes \vvm{\eta}{}{}  + \vvm{\hat{\Phi}}{}{} \cdot \vv{\hat{D}}{}{} : \vvm{\eta}{}{} \otimes \vvm{\chi}{}{} + \vvm{\hat{\Phi}}{}{} \cdot \vv{\hat{E}}{}{}: \vvm{\chi}{}{} \otimes \vvm{\chi}{}{},
 \label{eq:manApproxB}
 \end{align}
\end{subequations}
and eventually the reduced space of substructures displacements writes
\begin{equation}
\label{eq:manifGenerSubsA}
    \vv{d}{}{} = 
    \begin{bmatrix}
        \vv{u}{}{} \\ \vv{q}{}{}
    \end{bmatrix} 
    \approx 
    \begin{bmatrix}
        \vvm{\Phi}{}{}\vvm{\eta}{}{} + \vvm{\hat{\Phi}}{}{} \vvm{\Gamma}{\etahat}{}(\vvm{\eta}{}{},\vvm{\chi}{}{}) \\ \vvm{\Psi}{}{} \vvm{\chi}{}{}
    \end{bmatrix}.
\end{equation}
Upon defining the vector of substructure reduced coordinates as
\begin{equation}
        \vvm{\xi}{}{} = [\vvm{\eta}{}{ T},\vvm{\chi}{}{T}]^T \in \R{\ap{m}{}{}},
    \label{eq:decomposXi}
\end{equation}
the reduction manifold approximation is compactly written as
\begin{equation}
\label{eq:manifGenerSubs}
    \vv{d}{}{} 
    \approx 
    \vv{\Gamma}{d}{}(\vvm{\xi}{}{}) := \vv{L}{\Gamma}{}\vvm{\xi}{}{} + \vv{Q}{\Gamma}{}:(\vvm{\xi}{}{}\otimes \vvm{\xi}{}{}),
\end{equation}
where the linear manifold operator, $\vv{L}{\Gamma}{} \in \R{m\times m}$, writes
\begin{equation}
\label{eq:linearOpMan}
    \vv{L}{\Gamma}{} = 
    \begin{bmatrix}
        \vvm{\Phi}{}{} & \vvm{\hat{\Phi}}{}{} \vv{\hat{B}}{}{} \\ \zero & \vvm{\Psi}{}{}
    \end{bmatrix},
\end{equation}
and the quadratic operator, $\vv{Q}{\Gamma}{} \in \R{m\times m \times m}$, can be assembled from $\vvm{\hat{\Phi}}{}{} \cdot \vv{\hat{C}}{}{},\ \vvm{\hat{\Phi}}{}{} \cdot \vv{\hat{D}}{}{}, \ \text{and}\ \vvm{\hat{\Phi}}{}{} \cdot \vv{\hat{E}}{}{}$, as in Eq. \eqref{eq:manApproxB}.
\newl
It is interesting to note that if we truncate the manifold at linear level, we retrieve the reduction subspace used in classic CB reduction, as proved in \ref{sec:equivCB}.
This confirms the classic interpretation behind CB reduction, whereby high frequency FIMs are assumed to respond statically.
\subsubsection{Substructure Projected Equations}
\label{sec:projSubsEq}
The substructure ROM is obtained by Galerkin projection of the substructure equations onto the space tangent to the manifold.
Specifically, the dynamic residual vector  $\vvm{r}{}{}(\vvm{\xi}{}{},\vvmdot{\xi}{}{},\vvmddot{\xi}{}{}) \in \R{\ap{n}{d}{}}$ is defined by substituting Eq. \eqref{eq:manifGenerSubs} and its time derivatives into Eq. \eqref{eq:subsEq}, obtaining
\begin{equation}
    \label{eq:resSubs}
    \vvm{r}{}{}(\vvm{\xi}{}{},\vvmdot{\xi}{}{},\vvmddot{\xi}{}{}) :=  \vv{M}{}{} \vvmddot{\Gamma}{d}{}(\vvm{\xi}{}{},\vvmdot{\xi}{}{},\vvmddot{\xi}{}{}) + \vv{D}{}{} \vvmdot{\Gamma}{d}{}(\vvm{\xi}{}{},\vvmdot{\xi}{}{}) + \vv{K}{}{}\vvm{\Gamma}{d}{}(\vvm{\xi}{}{}) + \vv{f}{}{}(\vvm{\Gamma}{d}{}(\vvm{\xi}{}{})) - \vv{f}{ext}{}(t) - \vv{g}{}{},
\end{equation}
where 
\begin{subequations}
\label{eq:velAccMan}
    \begin{align}
            \vvmdot{\Gamma}{d}{}(\vvm{\xi}{}{},\vvmdot{\xi}{}{}) &= \vv{L}{\Gamma}{} \vvmdot{\xi}{}{} +  \vv{Q}{\Gamma,\tau}{}: \vvm{\xi}{}{} \otimes \vvmdot{\xi}{}{}, \\
            \vvmddot{\Gamma}{d}{}(\vvm{\xi}{}{},\vvmdot{\xi}{}{},\vvmddot{\xi}{}{}) &=  \vv{L}{\Gamma}{}\vvmddot{\xi}{}{} + \vv{Q}{\Gamma,\tau}{}: \vvmdot{\xi}{}{} \otimes \vvmdot{\xi}{}{} +
             \vv{Q}{\Gamma,\tau}{}: \vvm{\xi}{}{} \otimes \vvmddot{\xi}{}{},\\
            \vv{Q}{\Gamma,\tau}{} &= \vv{Q}{\Gamma}{} + \vv{Q}{\Gamma, 132}{}.
    \end{align}
\end{subequations}
Here,  tensor $\vv{Q}{\Gamma, 132}{}$ is obtained from $\vv{Q}{\Gamma}{}$ by permuting the elements in the second and third dimensions. 
Note that the first and second time derivatives of the manifold $\vvm{\Gamma}{d}{}$ include some convective terms that depend also on lower order time derivatives. 
For example, $ \vvmddot{\Gamma}{d}{}$ depends not only on $\vvmddot{\xi}{}{}$  but also on $\vvmdot{\xi}{}{}$ and $\vvm{\xi}{}{}$.
The convective terms arise from the quadratic component of the manifold, as also shown in \cite{jain2017}.\newl
Next, the residual is enforced to be orthogonal to the manifold tangent space, obtaining
\begin{equation}
    \pder{\vvm{\Gamma}{d}{T}}{\vvm{\xi}{}{}}(\vvm{\xi}{}{}) \cdot \vvm{r}{}{}(\vvm{\xi}{}{},\vvmdot{\xi}{}{},\vvmddot{\xi}{}{}) = \zero,
    \label{eq:galerProjSubs}
\end{equation}
where the projector is a linear function of $\vvm{\xi}{}{}$ and writes
\begin{equation}
    \label{eq:projector}
      \pder{\vvm{\Gamma}{d}{}}{\vvm{\xi}{}{}}(\vvm{\xi}{}{}) = \vv{L}{\Gamma}{} + \vv{P}{\Gamma}{}(\vvm{\xi}{}{}) =\vv{L}{\Gamma}{} + \vv{Q}{\Gamma,\tau}{} \cdot \vvm{\xi}{}{}.
\end{equation}
The reduced equations in \eqref{eq:galerProjSubs} are still not in the final form used for ROM time integration. 
In fact, a direct projection of the equations during time integration results in a computational cost that scales with the size of the FE model.
To solve this computational bottleneck, we perform an \textit{offline-online} decomposition, whereby the computationally intensive operations, required in the projection of HFM equations, are performed prior to time integration.
Specifically, we preassemble the ROM operators that multiply the reduced coordinates, as is usually done in Galerkin projection models that use a linear RB \cite{Saccani2025,Marconi2021,Mignolet2013}.
In that case, the reduced forces are cubic polynomials in the reduced coordinates, allowing for a tensorial representation which is extremely efficient to evaluate, and where the reduced nonlinear tensors are precomputed offline \cite{Mignolet2013}.
Here, instead, the projection on the quadratic manifold yields terms up to seventh order in $\vvm{\xi}{}{}$ for the elastic forces, as results from inserting the quadratic manifold expression into the cubic full elastic forces and projecting the equations on the linear projector in Eq. \eqref{eq:projector}.
Also in this case, the tensorial form could be adopted.
However, the higher polynomial order would undermine the efficient evaluation of the reduced forces, as well as the precomputation of the nonlinear tensors.
\newl
To improve efficiency of the ROM both in evaluation and model construction, we truncate the reduced nonlinear forces to third order. 
The validity of this assumption lies in the fact that we use the ROM to predict vibrations around an equilibrium point, so we expect  $\vvm{\xi}{}{}$ to be small.
As a result, we also expect that the relative importance of the monomial terms in $\vvm{\xi}{}{}$ in the reduced equation decreases with the monomial order. This  justifies truncation of higher order terms.
Additionally, also the quadratic manifold has local validity, as it is computed using perturbation theory.
For this reason, retaining all the monomial terms in the ROM equations is probably unnecessary, as higher order terms would become significant only when the magnitude of the reduced coordinates become large enough to invalidate the quadratic manifold hypothesis. Interestingly, a similar approximation is implicitly made in \cite{Fard2025}, where a reduction quadratic manifold is derived by assuming a quartic order potential function for ROM elastic forces.
\newl
The second simplifying assumption in the ROM derivation consists of neglecting the convective inertial terms in Eq. \eqref{eq:velAccMan}.
This hypothesis is consistent with neglecting the inertia and damping of high frequency FIMs in the manifold derivation, as detailed in Section \ref{sec:manDer}.
In fact, these convective terms are related to the inertia and damping of the quadratic component of the manifold, which has contributions only from fast FIMs, as can be seen from Eq. \eqref{eq:manApproxB}.
\newl
Eventually, the third and last simplifying assumption consists of neglecting the projection of the external load $\vv{f}{}{ext}(t)$ on the quadratic part of the manifold, $\vv{Q}{\Gamma}{}$. 
As for the previous point, this hypothesis is consistent with the manifold derivation, where we assumed that the forcing on high frequency FIMs is negligible (see Section \ref{sec:manDer}).
Again, this owes to the fact that the quadratic part of the manifold has contributions only from high frequency FIMs.
\newl
With the aforementioned simplifying assumptions, the projected residual equation in \eqref{eq:galerProjSubs} writes
\begin{equation}
    \vv{M}{r}{} \vvmddot{\xi}{}{} + \vv{D}{r}{} \vvmdot{\xi}{}{}  +\vv{K}{r}{} \vvm{\xi}{}{} + \vv{K}{r}{2} : (\vvm{\xi}{}{} \otimes \vvm{\xi}{}{}) + \vv{K}{r}{3} : (\vvm{\xi}{}{} \otimes \vvm{\xi}{}{} \otimes \vvm{\xi}{}{})  - \vv{f}{r}{ext}(\vvm{\xi}{}{},t) -\vv{g}{r}{} = \zero,
    \label{eq:subsROMeq}
\end{equation}
where the reduced linear operators write
\begin{equation}
        \vv{M}{r}{} = \vv{L}{\Gamma}{T} \vv{M}{}{} \vv{L}{\Gamma}{},\
        \vv{D}{r}{} = \vv{L}{\Gamma}{T} \vv{D}{}{} \vv{L}{\Gamma}{},\
        \vv{K}{r}{} = \vv{L}{\Gamma}{T} \vv{K}{}{} \vv{L}{\Gamma}{},\
         \vv{f}{r}{ext}(t) =  \vv{L}{\Gamma}{T} \vv{f}{ext}{}(t),\ \vv{g}{r}{} = \vv{L}{\Gamma}{T}\vv{g}{}{},
        \label{eq:redMassDampStiff}
\end{equation}
while the expressions for the reduced nonlinear stiffness tensors are
\begin{subequations}
\label{eq:tensorsExprROM}
    \begin{align}
         \vv{K}{r}{2} &= \vv{L}{\Gamma}{T} \cdot ( ( \vv{K}{}{2} \dotpr{21} \lman) \dotpr{21} \lman),\\
         \vv{K}{r}{3} &= \vv{L}{\Gamma}{T} \cdot (( \vv{K}{}{2} \dotpr{21} \lman) \dotpr{21} \qman ) + \vv{L}{\Gamma}{T} \cdot (( \vv{K}{}{2} \dotpr{21} \qman) \dotpr{21} \lman ) + \vv{L}{\Gamma}{T} \cdot ((( \vv{K}{}{3} \dotpr{21} \lman) \dotpr{21} \lman ) \dotpr{21} \lman).
    \end{align}
\end{subequations}
We refer the reader to \ref{sec:redNlTens} for a formal derivation of the nonlinear ROM operators.
Notably, Eq. \eqref{eq:subsROMeq} has the same structure of the ROM equations obtained by projecting the nonlinear FE equations on a linear RB, guaranteeing optimal efficiency.
\subsubsection{On the Lagrangian Properties of the Substructure ROM}
In this section, we show that the substructure ROM equations in \eqref{eq:subsROMeq} inherit the \textit{Lagrangian Structure} \cite{Farhat2015} underlying the starting FE model equations.
In particular, we show that (i) the reduced mass, damping, and stiffness matrices are positive definite, and that (ii) the reduced elastic forces can be derived from a reduced potential function.
If these conditions are satisfied, reduced kinetic and potential energies, as well as a dissipation function, can be defined.
Preservation of the Lagrangian structure is a sought-after feature when devising new model reduction techniques, as it results in ROMs whose time-integration is numerically stable for a proper choice of time-integration schemes \cite{Farhat2015}.
The ROMs constructed using the Galerkin Projection on a linear RB are known to be structure preserving \cite{Farhat2015}. 
In this work, however, we obtain the substructure ROM equations by Galerkin projection onto the tangent space of a quadratic manifold and by applying the simplifying assumptions detailed in Section \ref{sec:projSubsEq}.
This procedure yields a structure preserving ROM, as shown in the following.
\newl
The  reduced mass, damping, and stiffness matrices are positive definite, as they are obtained by left and right projection of their positive definite counterpart on the linear manifold operator $\vv{L}{\Gamma}{}$, as shown in Eq. \eqref{eq:redMassDampStiff}.
Thus, the ROM kinetic energy and dissipation function, $\mathcal{T}$ and $\mathcal{D}$, can be defined as 
\begin{subequations}
\begin{align}
    \mathcal{T}_{r}(\vvmdot{\xi}{}{}) = &\frac{1}{2}\vvmdot{\xi}{}{T}\vv{M}{r}{}\vvmdot{\xi}{}{}, \\
     \mathcal{D}_{r}(\vvmdot{\xi}{}{}) = & \frac{1}{2}\vvmdot{\xi}{}{T}\vv{D}{r}{}\vvmdot{\xi}{}{}. 
\end{align}
\end{subequations}
Moreover, the substructure ROM admits an elastic reduced potential energy of the form
\begin{equation}
    \mathcal{V}_r(\vvm{\xi}{}{}) = \frac{1}{2} \vvm{\xi}{}{} \cdot \vv{K}{r}{} \cdot \vvm{\xi}{}{} + \frac{1}{3}\vvm{\xi}{}{} \cdot \vv{K}{r}{2} : (\vvm{\xi}{}{} \otimes \vvm{\xi}{}{}) +  \frac{1}{4}\vvm{\xi}{}{} \cdot \vv{K}{r}{3} \cdddot   (\vvm{\xi}{}{} \otimes \vvm{\xi}{}{}  \otimes \vvm{\xi}{}{}), 
\end{equation}
where $\vv{K}{r}{},\vv{K}{r}{2}$, and $\vv{K}{r}{3}$ are the substructure reduced stiffness matrix and nonlinear stiffness tensors reported in Eq. \eqref{eq:tensorsExprROM}.
We refer the reader to \ref{sec:redPotEn} for a proof of this statement.
\newl
As a consequence,  along a solution of the ROM equation, $\vvm{\xi}{}{}(t)$, the following power balance holds
\begin{equation}
    \frac{d}{dt} \mathcal{E}_r(\vvm{\xi}{}{}(t),\vvmdot{\xi}{}{}(t)) = -\vvmdot{\xi}{}{}(t)^T \pder{}{\vvmdot{\xi}{}{}} \mathcal{D}_r(\vvmdot{\xi}{}{}(t)) + \vvmdot{\xi}{}{}(t)^T\vv{f}{}{r}(t) + \vvmdot{\xi}{}{}(t)^T \vv{g}{r}{}(t),
\end{equation}
where the total energy of the reduced system $\mathcal{E}_r(\vvm{\xi}{}{}(t),\vvmdot{\xi}{}{}(t))$ is defined as
\begin{equation}
    \mathcal{E}_r(\vvm{\xi}{}{},\vvmdot{\xi}{}{}) = \mathcal{T}_{r}(\vvmdot{\xi}{}{}) + \mathcal{V}_r(\vvm{\xi}{}{}).
\end{equation}
Note the presence of the energy related to the interface forces.
This term cancels out after assembling the substructures ROM, as the interface forces are pairs of action-reaction forces and compatibility of displacements is strongly satisfied at the interface.
\subsection{Computation of the Substructure ROM}
\subsubsection{Computation of the Manifold}
\label{sec:manComp}
As presented in Section \ref{sec:manDer} and detailed in \ref{app:secOrdS}, the substructure manifold is derived using perturbation theory.
The resulting equations for the Taylor series coefficients vectors, stored in $\vv{\hat{B}}{}{}$, $\vv{\hat{C}}{}{}$, $\vv{\hat{D}}{}{}$, and $\vv{\hat{E}}{}{}$ (see Eq. \eqref{eq:manFast}), take the general form
\begin{equation}
    \vvm{\hat{\Phi}}{}{T} \vv{K}{uu}{} \vvm{\hat{\Phi}}{}{} \vvm{\hat{\eta}}{}{*} = \vvm{\hat{\Phi}}{}{T} \vvm{f}{}{*},
    \label{eq:genericEqCoeff}
\end{equation} 
where $\vvm{\etahat}{}{*} \in \R{\ap{n}{\hat{\phi}}{}}$ and $\vvm{f}{}{*} \in \R{\ap{n}{u}{}}$ are, respectively, a generic manifold component vector and right hand side forcing arising from perturbation.
In other words, $\vvm{\eta}{}{*}$ can be any of $\vv{\hat{B}}{}{}(:,i)$, $\vv{\hat{C}}{}{}(:,i,j)$, $\vv{\hat{D}}{}{}(:,i,j)$, and $\vv{\hat{E}}{}{}(:,i,j)$, where we use Matlab notation to extract the component vectors.
The direct solution of Eq.\eqref{eq:genericEqCoeff} is not a feasible option, as this would imply the computation of all the high frequency FIMs $\vvm{\hat{\Phi}}{}{}$, which is intractable for large FE models. \newl
A workaround to this problem is found by exploiting stiffness orthogonality of FIMs. 
In fact, $\vvm{\hat{\Phi}}{}{} \vvm{\hat{\eta}}{}{*}$ that satisfies Eq. \eqref{eq:genericEqCoeff} can be computed as
\begin{equation}
    \vvm{\hat{\Phi}}{}{}\vvm{\hat{\eta}}{}{*} = (\vv{I}{}{} - \vvm{\Phi}{}{}\vvm{\Phi}{}{T}\vv{M}{uu}{})\vv{K}{uu}{-1}\vvm{f}{}{*}.
    \label{eq:exprr}
\end{equation}
Upon substitution of Eq. \eqref{eq:exprr} into Eq. \eqref{eq:genericEqCoeff} and using stiffness orthogonality of FIMs, it is easy to see that this is the solution to Eq. \eqref{eq:genericEqCoeff}. 
We remark that the computation of the different manifold coefficients requires the inversion of the same linear operator $\vv{K}{uu}{}$ , thus its
upfront factorization can yield computational speed-ups.
In this case, the Cholesky factorization can be used, as the stiffness matrix is symmetric positive definite.
The reader should note that this alternative computational procedure does not require the computation of all FIMs, but instead of only low frequency FIMs, which are usually few.
Note that by solving Eq. \eqref{eq:exprr}, we do not compute $\vvm{\eta}{}{*}$, but only $\vvm{\hat{\Phi}}{}{} \vvm{\hat{\eta}}{}{*}$.
However, this is not an issue because the manifold operators $\vv{\hat{B}}{}{}$, $\vv{\hat{C}}{}{}$, $\vv{\hat{D}}{}{}$, and $\vv{\hat{E}}{}{}$ are eventually premultiplied by $\vvm{\hat{\Phi}}{}{}$, as shown in Eq. \eqref{eq:manApproxB}.
\newl
Another procedure that requires particular care in the computation of the manifold is the assembly of the forcing vector $\vvm{f}{}{*}$ on the right hand side of Eqs. \eqref{eq:genericEqCoeff}. 
To compute the coefficients of the manifold at leading order, $\vvm{f}{}{*}$ can be constructed by projecting the full order stiffness matrix $\vv{K}{uu}{}$ on the basis of low frequency FIMs, $\vvm{\Phi}{}{}$, and of SMs, $\vvm{\Psi}{}{}$. 
Conversely, a different procedure is followed in the computation of the second order coefficients.
In this case, as can be seen from \ref{app:secOrdS}, vector $\vvm{f}{}{*}$ takes the general form
\begin{align}
    \vvm{f}{}{*} &= \sum_i \vv{K}{uuu}{2} : (\vv{v}{u}{i} \otimes \vv{w}{u}{i} +  \vv{w}{u}{i} \otimes \vv{v}{u}{i} ) + \nonumber \\ 
     &+ \sum_j \vv{K}{uqq}{2} : (\vv{v}{q}{j} \otimes \vv{w}{q}{j} + \vv{w}{q}{j}  \otimes \vv{v}{q}{j} ) + \nonumber \\
    &+ \sum_k \vv{K}{uuq}{2} : \vv{v}{u}{k} \otimes \vv{v}{q}{k} .
    \label{eq:rhsSecOrd}
\end{align}
Here, the tensors associated to second order nonlinear forces are contracted two times on vectors of displacements partitioned in interface and internal DoFs, which we denote with $\vv{v}{}{}$ and $\vv{w}{}{}$.
The computation of vector $\vvm{f}{}{*}$ by performing tensor-vector contractions is computationally infeasible, as it requires the assembly of the full order tensors, $\vv{K}{uuu}{2},\vv{K}{uqq}{2}$, and $\vv{K}{uuq}{2}$. This operation is extremely memory consuming even for moderately small FE models.\newl
To overcome this limitation, we assemble $\vvm{f}{}{*}$ using the tangent stiffness derivative approach, following a similar procedure to the one used in \cite{KaramoozMahdiabadi2021,jain2017,Marconi2021} for the computation of the SMDs.
Specifically, we consider the following relations
\begin{subequations}
\begin{align}
\label{eq:tangStiffDer}
        \vv{K}{uuu}{2} : (\vv{v}{u}{} \otimes \vv{w}{u}{} +\vv{w}{u}{} \otimes \vv{v}{u}{}) &= \frac{\partial}{\partial \epsilon} \frac{\partial \vv{f}{u}{el}}{\partial \vv{u}{}{}}(\epsilon \vv{v}{u}{},\zero) \bigg{|}_{\epsilon = 0} \cdot \vv{w}{u}{},\\  
         \vv{K}{uqq}{2} : (\vv{v}{q}{} \otimes \vv{w}{q}{} + \vv{w}{q}{} \otimes \vv{v}{q}{}) &= \frac{\partial}{\partial \epsilon} \frac{\partial \vv{f}{u}{el}}{\partial \vv{q}{}{}}(\zero, \epsilon \vv{v}{q}{}) \bigg{|}_{\epsilon = 0} \cdot \vv{w}{q}{},\\
        \vv{K}{uuq}{2} : \vv{v}{u}{} \otimes \vv{v}{q}{} &= \frac{\partial}{\partial \epsilon} \frac{\partial \vv{f}{u}{el}}{\partial \vv{q}{}{}}(\epsilon \vv{v}{u}{}, \zero) \bigg{|}_{\epsilon = 0} \cdot \vv{v}{q}{} = \frac{\partial}{\partial \epsilon} \frac{\partial \vv{f}{u}{el}}{\partial \vv{u}{}{}}(\zero,\epsilon \vv{v}{q}{}) \bigg{|}_{\epsilon = 0} \cdot \vv{v}{u}{},
\end{align}
\end{subequations}
where
\begin{equation}
    \vv{f}{u}{el}(\vv{u}{}{},\vv{q}{}{}) = \vv{K}{uu}{}\vv{u}{}{} + \vv{K}{uq}{} \vv{q}{}{} + \vv{f}{u}{}(\vv{u}{}{},\vv{q}{}{})
\end{equation}
are the substructure elastic forces on internal DoFs and
$\frac{\partial \vv{f}{u}{el}}{\partial \vv{u}{}{}}(\vv{u}{}{}, \vv{q}{}{}) \text{ and } \frac{\partial \vv{f}{u}{el}}{\partial \vv{q}{}{}}(\vv{u}{}{}, \vv{q}{}{}) $ are two partitions of the  tangent stiffness matrix.
We refer the reader to
\ref{sec:tangStiffRel} for a proof of these expressions.
In this setting, the tangent stiffness derivatives in Eq. \eqref{eq:tangStiffDer} can be computed with a specifically written code or using finite differences. 
In this second case, for example, we can write
\begin{equation}
    \frac{\partial}{\partial \epsilon} \frac{\partial \vv{f}{u}{el}}{\partial \vv{u}{}{}}(\epsilon \vv{v}{u}{}, \zero) \bigg{|}_{\epsilon = 0} \approx 
   \left( \frac{\partial \vv{f}{u}{el}}{\partial \vv{u}{}{}}(+h \vv{v}{u}{}, \zero)  -  \frac{\partial \vv{f}{u}{el}}{\partial \vv{u}{}{}}(-h \vv{v}{u}{}, \zero) \right)(2h)^{-1},
\end{equation}
where a first order central finite difference scheme with incremental step $h$ is used.\newl
The procedure for the computation of the manifold used to reduce the substructures is summarized in Algorithm \ref{alg:ROMconstr}.
\begin{algorithm}[H]
\caption{Substructure Manifold Computation}
\textbf{Input:} \texttt{SubsFeModel}\sscr{a}, $n_{\phi}$\sscr{b}, $\vvm{\Psi}{}{}$\sscr{c} \\ 
\textbf{Output}: $\vv{L}{\Gamma}{}, \vv{Q}{\Gamma}{}$  .
\label{alg:ROMconstr}
\begin{algorithmic}[1] 
    \Statex \textit{Computation of FIMs} 
    \State $\vv{K}{uu}{} \gets  \texttt{SubsFeModel.Kuu}$ \sscr{d}
    \State $\vv{M}{uu}{} \gets  \texttt{SubsFeModel.Muu}$  \sscr{e}
    \State $\vvm{\Phi}{}{} \gets  \texttt{eigs}(\vv{M}{uu}{},\vv{K}{uu}{},n_{\phi})$ \sscr{f}
    \Statex
    \Statex \textit{Computation of the Linear Manifold Component}
    \State $\vv{K}{\eta\eta}{} \gets \vvm{\Phi}{}{T} \vv{K}{uu}{} \vvm{\Phi}{}{} $
    \State $\vvm{\hat{\Phi}}{}{} \vv{\hat{B}}{}{}\gets   \vv{K}{uu}{-1}\vv{K}{uq}{}\vvm{\Psi}{}{} - \vvm{\Phi}{}{} \vv{K}{\eta \eta}{-1} \vv{K}{uq}{}\vvm{\Psi}{}{} $
    \State $\vv{L}{\Gamma}{} \gets [\vvm{\Phi}{}{},\vvm{\hat{\Phi}}{}{} \vv{\hat{B}}{}{}; \zero, \vvm{\Psi}{}{}]$
    \Statex
    \Statex \textit{Computation of the Quadratic Manifold Component}
    \State $\vv{F}{}{*} \gets \texttt{RHS2ndOrdCoeff}(\texttt{SubsFEModel},\vvm{\Phi}{}{},\vvm{\hat{\Phi}}{}{}\vv{\hat{B}}{}{},\vvm{\Psi}{}{})$ \sscr{g}
    \Statex
    \State $\vv{X}{}{*} \gets \vv{K}{uu}{}\backslash\vv{F}{}{*}$ \sscr{h}
    \State $\vv{X}{}{*} \gets (\vv{I}{}{}-\vvm{\Phi}{}{}\vvm{\Phi}{}{T}\vv{M}{uu}{})\vv{X}{}{*}$
    \Statex
    \State $\vv{Q}{\Gamma}{} \gets \texttt{assembleQuadMan}(\vv{X}{}{*})$\sscr{i}
\end{algorithmic}
\end{algorithm}
\vspace{-10pt}
\footnotesize
\noindent 
\sscr{a}\texttt{SubsFeModel} is the substructure FE model.\\
\sscr{b}Number of FIMs used in substructure reduction.\\
\sscr{c}RB for interface DoFs.\\
\sscr{d} \texttt{SubsFeModel.Kuu} returns the internal DoFs stiffness matrix.\\
\sscr{e} \texttt{SubsFeModel.Muu} returns the internal DoFs mass matrix.\\
\sscr{f} Function \texttt{eigs} returns the first $n_\phi$ VMs.\\
\sscr{g} Function \texttt{RHS2ndOrdCoeff}  computes, following Algorithm \ref{alg:ROMconstrRHS} in  \ref{sec:appAlg}, the right hand side forces needed in the computation of the quadratic  manifold coefficients, as defined by Eq.\eqref{eq:rhsSecOrd}. \\
\sscr{h} Invert the linear system $\vv{K}{uu}{}\vv{x}{}{*}=\vv{F}{}{*}$.\\
\sscr{i} Function \texttt{assembleQuadMan} reshapes the quadratic manifold components in $\vv{X}{}{*}$ in tensor $\vv{Q}{}{\Gamma}$.
\normalsize
\subsubsection{Assembly of Substructure ROM Tensors using Element Level Projection}
The assembly of the nonlinear substructure ROM tensors in Eq. \eqref{eq:subsROMeq} is performed using an element level projection strategy. 
A direct assembly of the substructure full order tensors, $\vv{K}{2}{}$ and $\vv{K}{3}{}$, and their subsequent contraction with the manifold operators using the formulas in Eq. \eqref{eq:tensorsExprROM} results in large computational memory consumption.
In the element level projection method, however, the assembly of full order tensors is avoided.
Instead, only the element level tensors are computed and directly contracted with the manifold operators restricted to the element variables.
Specifically, for a generic element $e$ in the element set $\mathcal{E}$ of a generic substructure, we can define the vector of nodal variables, $\vv{d}{e}{} \in \R{n_e}$, and the element localization matrix $\vvm{L}{e}{}\in \R{n_e\times n}$, such that $\vv{d}{e}{} = \vvm{L}{e}{} \vv{d}{}{}$. 
Next, we consider the manifold operators restricted to element variables
\begin{equation}
    \vv{L}{\Gamma,e}{} = \vvm{L}{e}{} \vv{L}{\Gamma}{} \in \R{n_e \times m},\quad 
    \vv{Q}{\Gamma,e}{} = \vvm{L}{e}{} \cdot \vv{Q}{\Gamma}{} \in \R{n_e \times m \times m},
\end{equation}
and the element stiffness matrix and nonlinear tensors,
$\vv{K}{e}{} \in \R{n_e \times n_e}$, $\vv{K}{e}{2} \in \R{n_e \times n_e \times n_e}$, and $\vv{K}{e}{3} \in \R{n_e \times n_e \times n_e \times n_e}$.
In the element level projection scheme, the nonlinear tensors are computed using the following formulas, derived from Eq. \eqref{eq:tensorsExprROM}:
\begin{subequations}
\label{eq:tensorsExprROMpr}
    \begin{align}
         \vv{K}{r}{2} &= \sum_{e \in \mathcal{E} }^{N_e} \vv{L}{\Gamma,e}{T} \cdot ( ( \vv{K}{e}{2} \dotpr{21} \lmane) \dotpr{21} \lmane),\\
         \vv{K}{r}{3} &=  \sum_{e \in \mathcal{E}} \vv{L}{\Gamma,e}{T} \cdot (( \vv{K}{e}{2} \dotpr{21} \lmane) \dotpr{21} \qmane ) + \vv{L}{\Gamma,e}{T} \cdot (( \vv{K}{e}{2} \dotpr{21} \qmane) \dotpr{21} \lmane ) + \\ &\quad \quad + \vv{L}{\Gamma,e}{T} \cdot ((( \vv{K}{e}{3} \dotpr{21} \lmane) \dotpr{21} \lmane ) \dotpr{21} \lmane). \nonumber
    \end{align}
\end{subequations}
The element level projection approach relies on specialized FE programs that can return the nonlinear tensors related to nonlinear elastic forces. 
Usually, this tensorial decomposition is not an output of commercial FE software.
To overcome this limitation, the method presented here could be potentially integrated with non-intrusive tensor identification methods in the same spirit of the Enforced Displacements Method or the Enhanced Enforced Displacements Method \cite{Mignolet2013}.
However, such algorithms would need to be adapted to allow the identification of polynomial approximate reduced forces and to accommodate the projection on quadratic manifold components.
\subsection{ROM Assembly}
\label{sec:romass}
The final form of the ROM is obtained by assembling the substructure reduced equations in \eqref{eq:subsROMeq}.
Since these equations are substructure related, we report them here reinstating the superscript '$(s)$':
\begin{equation}
    \vv{M}{r}{(s)} \vvmddot{\xi}{}{(s)} + \vv{D}{r}{(s)} \vvmdot{\xi}{}{(s)}  +\vv{K}{r}{(s)} \vvm{\xi}{}{(s)} + \vv{K}{r}{2(s)} : (\vvm{\xi}{}{(s)} \otimes \vvm{\xi}{}{(s)}) + \vv{K}{r}{3(s)} : (\vvm{\xi}{}{(s)} \otimes \vvm{\xi}{}{(s)} \otimes \vvm{\xi}{}{(s)})  - \vv{f}{r}{ext(s)}(\vvm{\xi}{}{(s)},t) -\vv{g}{r}{(s)} = \zero,
    \label{eq:subsROMeq_s}
\end{equation}
where $\vvm{\xi}{}{(s)} \in \R{m^{(s)}}$ is the vector of substructure reduced variables.
A vector of global reduced coordinates, $\vvm{\tilde{\xi}}{}{} \in \R{m}$, is defined from the vectors of substructure reduced variables by assigning the same variable to interface displacements that are shared across substructures.
If interface reduction is used, we assume that the interface RBs are compatible.
The substructure reduced variables could be extracted using localization matrices from the vector of global reduced coordinates $\vvm{\tilde{\xi}}{}{}$ as
\begin{equation}
    \vvm{\xi}{}{(s)} = \vv{L}{r}{(s)} \vvm{\tilde{\xi}}{}{},
    \label{eq:localSubs}
\end{equation}
where $\vv{L}{r}{(s)} \in \R{m^{(s)} \times m}$.
The ROM equations are obtained by substituting Eq. \eqref{eq:localSubs} in Eq. \eqref{eq:subsROMeq_s}, projecting the substructure equations on the localization matrices $\vv{L}{r}{(s)}$, and by summing over the different substructures composing the assembly.
Thus, the final form of the ROM is
\begin{equation}
    \vv{M}{r}{} \ddot{\bm{\tilde{\xi}}} + \vv{D}{r}{} \vvmdot{\tilde{\xi}}{}{}  +\vv{K}{r}{} \vvm{\tilde{\xi}}{}{} + \vv{K}{r}{2} : (\vvm{\tilde{\xi}}{}{} \otimes \vvm{\tilde{\xi}}{}{}) + \vv{K}{r}{3} : (\vvm{\tilde{\xi}}{}{} \otimes \vvm{\tilde{\xi}}{}{} \otimes \vvm{\tilde{\xi}}{}{})  - \vv{f}{r}{ext}(\vvm{\tilde{\xi}}{}{},t) = \zero,
    \label{eq:subsROMfinal}
\end{equation}
where
\begin{subequations}
    \begin{align}
        &\vv{M}{r}{} = \sum_{s=1}^{N_s} \vv{L}{r}{(s)T} \vv{M}{r}{(s)}  \vv{L}{r}{(s)}, \ \  \vv{D}{r}{} = \sum_{s=1}^{N_s} \vv{L}{r}{(s)T} \vv{D}{r}{(s)}  \vv{L}{r}{(s)}, \ \ \vv{K}{r}{} = \sum_{s=1}^{N_s} \vv{L}{r}{(s)T} \vv{K}{r}{(s)}  \vv{L}{r}{(s)}, \\
        &\vv{K}{r}{2} =  \sum_{s=1}^{N_s} \vv{L}{r}{(s)T} \cdot ((\vv{K}{r}{2(s)}\cdot_{21} \vv{L}{r}{(s)} )\cdot_{21} \vv{L}{r}{(s)} ), \\
        &\vv{K}{r}{3} =  \sum_{s=1}^{N_s} \vv{L}{r}{(s)T} \cdot (((\vv{K}{r}{3(s)}\cdot_{21} \vv{L}{r}{(s)} )\cdot_{21} \vv{L}{r}{(s)})\cdot_{21} \vv{L}{r}{(s)} ),\\
        &\vv{f}{r}{ext}(t) =  \sum_{s=1}^{N_s} \vv{L}{r}{(s)T} \vv{f}{r}{ext(s)}(t) 
        \label{eq:redForc11}
    \end{align}
\end{subequations}
The reader should note that the ROM assembly procedure eliminates the reaction forces at the interface between different substructures:
\begin{equation}
    \sum_{i=1}^{N_s} \vv{L}{r}{(s)T} \vv{g}{r}{(s)} = 0.
\end{equation}
This is because we assume that the interfaces of different connected substructures are reduced with compatible interface modes.
\newl
Time integration of Eq. \eqref{eq:subsROMfinal} entails Newton-Raphson iterations, which rely on the computation of its Jacobian with respect to the reduced variables.
The expression for the ROM Jacobian writes
\begin{equation}
  \pder{\vvm{r}{}{}}{\vvm{\tilde{\xi}}{}{}} = \vv{K}{r}{} + \vv{K}{r}{2\text{t}} \cdot \vvm{\tilde{\xi}}{}{} + \vv{K}{r}{3\text{t}} : ( \vvm{\tilde{\xi}}{}{} \otimes  \vvm{\tilde{\xi}}{}{} ) , \quad   \pder{\vvm{r}{}{}}{\vvmdot{\tilde{\xi}}{}{}} = \vv{D}{r}{},  \quad   \pder{\vvm{r}{}{}}{\vvmddot{\tilde{\xi}}{}{}} = \vv{M}{r}{}.
\end{equation}
where
\begin{equation}
    \vv{K}{r}{2\text{t}} = \vv{K}{r}{2} +\vv{K}{r, 132}{2}, \quad  \vv{K}{r}{3\text{t}} = \vv{K}{r}{3} + \vv{K}{r, 1324}{3} +\vv{K}{r, 1243}{3}.
\end{equation}
In this way, the ROM equations can be time-integrated using the classic \textit{Newmark-$\beta$} time integration scheme \cite{geradin2015mechanical}.
\section{Applications} \label{sec:applications}
To test the validity of the proposed method, we performed four different numerical experiments, benchmarking the nonlinear substructure ROM with full order FE simulations.
The numerical examples considered are flat and curved von K\'arm\'an beams, a flat panel, and a MEMS gyroscope.
The FE implementation used in all numerical experiments is based on \textit{YetAnotherFEcode} \cite{yafec}, our in-house written Matlab code.
\subsection{Flat von K\'arm\'an Beam}
\label{flatbeam}
\begin{figure}[h]
    \centering
    \includegraphics[width=1\linewidth]{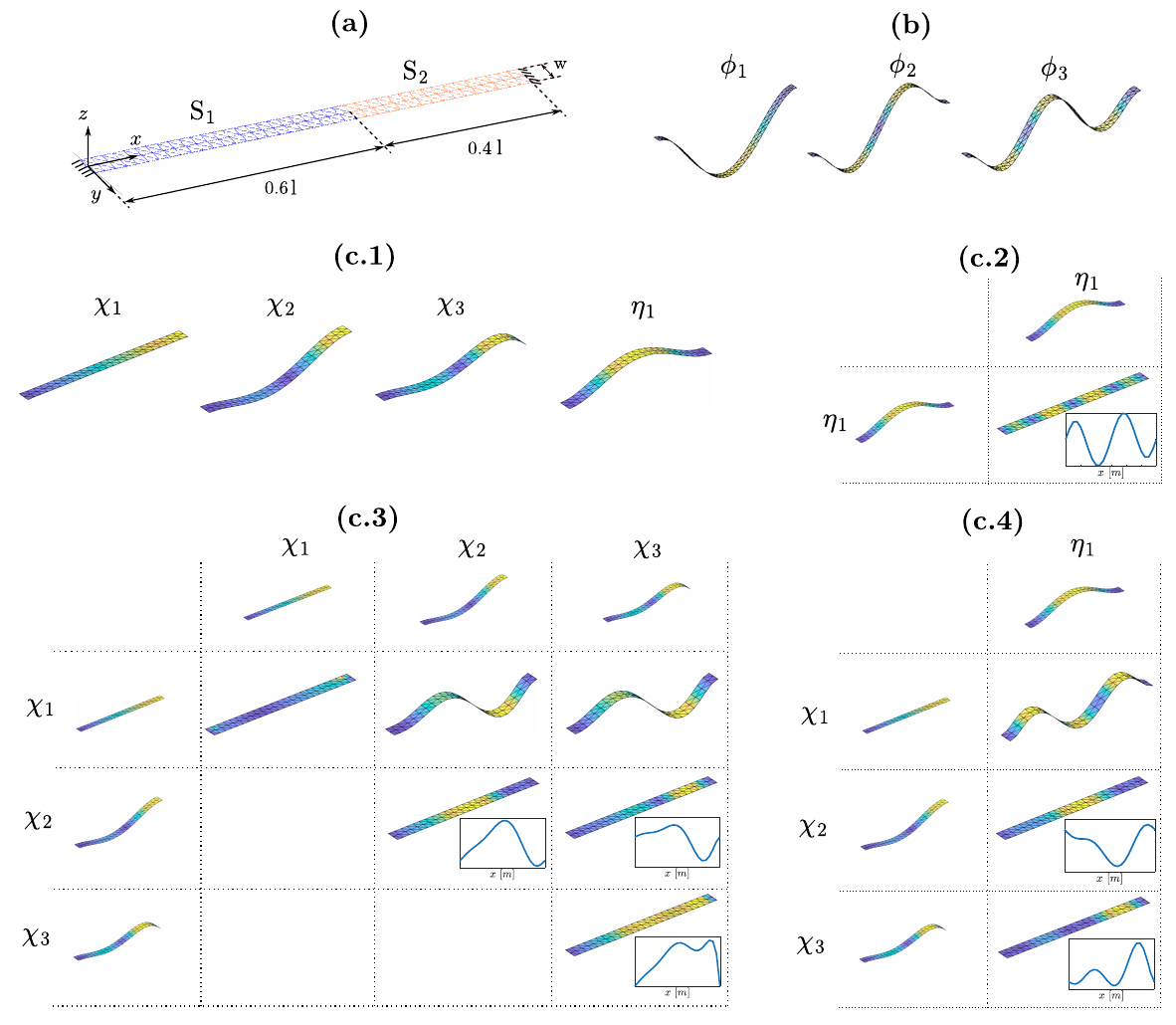}
    \caption{Flat clamped-clamped beam withvon K\'arm\'an nonlinearities. In (a), FE assembly partitioned in two substructures $S_1$ and $S_2$. In (b), the first three VMs of the assembled structures. In (c), substructure manifold components for substructure $S_1$. In (c.1) we show the linear terms, whereas, in (c.2) the quadratic terms in $\eta$, in (c.3) the quadratic terms in $\chi$, and in (c.4) the quadratic cross-coupling terms in $\eta$ and $\chi$. }
    \label{fig:flatBeamModel}
\end{figure}

\begin{figure}[h]
    \centering
    \includegraphics[width=1\linewidth]{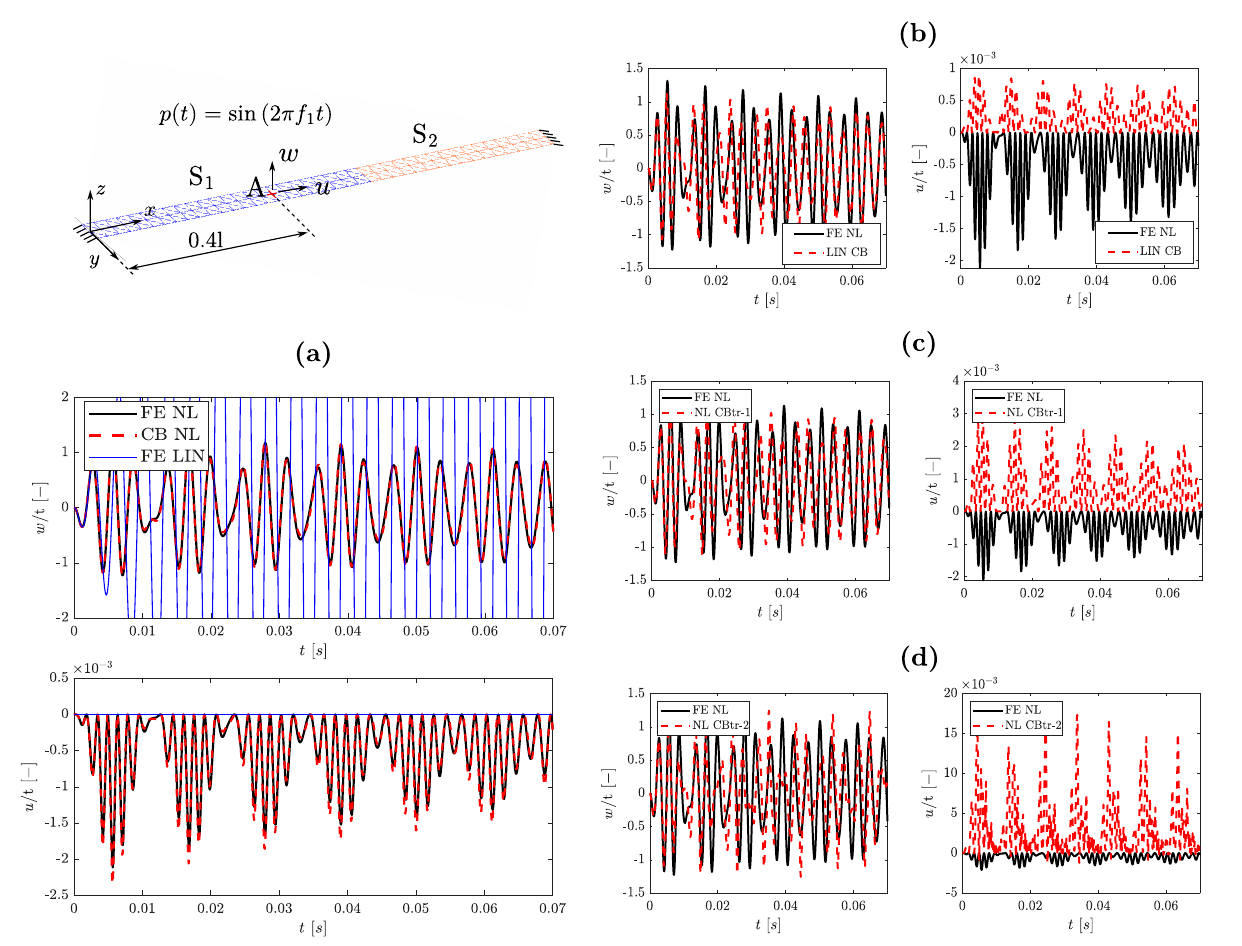}
    \caption{Accuracy of the CB ROM in predicting the vibration of the flat clamped-clamped von K\'arm\'an beam subjected to uniform pressure loading, varying harmonically at the natural frequency. In (a), out-of-plane (top) and in-plane (bottom) response at point $\text{A}$, computed using the full FE nonlinear model (FE-NL), the NL-CB ROM (CB NL), and the linearized model (FE-LIN). In (b), (c), and (d), ROM approximation compared to the full model (FE-NL), when the quadratic manifold part is set to zero (b), the quadratic terms in $\chi$ are set to zero (c), and the quadratic cross-coupling terms in $\chi$ and $\eta$ are neglected.}
    \label{fig:flatBeamResp}
\end{figure}

The first investigated example is a flat clamped-clamped beam with geometric nonlinearities stemming from the von K\'arm\'an strain measure.
The beam has length $\text{l} = 100\ \text{mm}$, width $\text{w} = 5\ \text{mm} $, and a rectangular cross section of thickness $\text{t} = 0.5 \ \text{mm}$.  Linear elastic material with an elastic modulus $E= 210\ \text{GPa}$, Poisson ratio $\nu = 0.33$, and density $\rho = 7800 \ \text{Kg}/\text{m}^3$,
is considered.
The FE discretization was performed using $160$ linear shell elements, each with $18$ DoFs, yielding the FE mesh reported in Fig.  \ref{fig:flatBeamModel}a. To prevent cross section torsion, the rotations about the $y$ axes of the nodes along $y= 0$ were constrained to zero. After applying boundary conditions, the total number of DoFs is $585$.
The first three VMs of the structures are reported in Fig. \ref{fig:flatBeamModel}b, and their associated frequencies are $f_1=269.5\ \text{Hz},\ f_2=742.8\ \text{Hz}, \text{and}\  f_3=1456.8\ \text{Hz}$.
Rayleigh damping with $\alpha = 24.85$ and $\beta = 3.15\cdot 10^{-6} $ was considered. 
These values were obtained by enforcing a damping ratio of $1\%$ to the first two structural modes.
\newl
As shown in Fig. \ref{fig:flatBeamModel}a, the FE assembly was divided into two substructures, $S_1$ and $S_2$, by cutting the beam at $0.6\ \text{l}$.
Then, the NL-CB ROM was computed using  one FIM for substructure.
The interface was reduced leveraging the physics of the problem, assuming rigid interface motion, as expected in beams.
This interface reduction technique is also known as the Virtual Node method \cite{Krattiger2019}.
Specifically, three interface modes were used, two of which correspond to uniform displacements of the interface along $z$ and $x$ axes, whereas the third one is a uniform rotation of the interface about the $y$ axis.
This results in an assembled substructure ROM of dimension five.
The linear and nonlinear components of the substructure  manifold of substructure $S_1$ are plotted in Fig. \ref{fig:flatBeamModel}c.
By looking at the nonlinear manifold components in Fig. \ref{fig:flatBeamModel}c.2 and Fig. \ref{fig:flatBeamModel}c.3, it is interesting to note that the quadratic coupling terms between out of plane modes induce second order in plane displacements  - see for example the coupling between the amplitude pairs $\chi_2-\chi_2$, $\chi_2-\chi_3$, $\chi_3-\chi_3$,  $\chi_2-\eta_1$,  $\chi_3-\eta_1$, and $\eta_1-\eta_1$.
This coupling mechanism, known as  \textit{bending-stretching} coupling,  has also been reported in other numerous studies \cite{Mignolet2013, KaramoozMahdiabadi2021, Kim2013}.
\newl
The NL-CB ROM accurately predicts the first two natural frequencies of the flat beam, respectively identified with $0.42 \%$ and $0.40 \%$ frequency error.
To test the accuracy of the ROM in the nonlinear regime, we computed the response of the beam to uniform pressure, oscillating harmonically at the first natural frequency as $p(t) = 20 \sin{(2\pi f_1 t)}\ \text{kPa}$.
We time integrate the ROM, the full nonlinear and the linearized FE models, comparing the displacement time history at node $A$ in Fig. \ref{fig:flatBeamResp}.
Time integration of the three models was performed with the Newmark-$\beta$ time integration scheme \cite{geradin2015mechanical}, using a constant time step of $7.42\times 10^{-5} \text{s}$. The NL-CB ROM is in good agreement with the full FE model and accurately approximates out-of-plane and in-plane displacements, as shown in Fig.  \ref{fig:flatBeamResp}a.
The response is remarkably nonlinear, as can be seen from the comparison with the linear model. 
Additionally, to provide further insight on the method, we time-integrated other three different substructure ROMs, obtained from the NL-CB ROM by setting to zero (i) the quadratic part, (ii) the quadratic part multiplying the amplitudes of SMs only, and (iii) the quadratic coupling terms between FIMs and SMs.
These approximations are shown respectively in Figs.  \ref{fig:flatBeamResp}b, \ref{fig:flatBeamResp}c, and \ref{fig:flatBeamResp}d.
As can be seen, all the quadratic terms are essential to capture the nonlinear behavior of the beam, especially the in-plane motion. 
Without the quadratic part (case i), the substructure ROM is equivalent to the ROM constructed by Galerkin projection of the nonlinear equations using the conventional CB RB.
In this case, the ROM predicts an overly stiff response, as the bending-stretching coupling is not well represented by the CB RB.
\subsection{Curved von K\'arm\'an Beam}
\begin{figure}[h]
    \centering
    \includegraphics[width=1\linewidth]{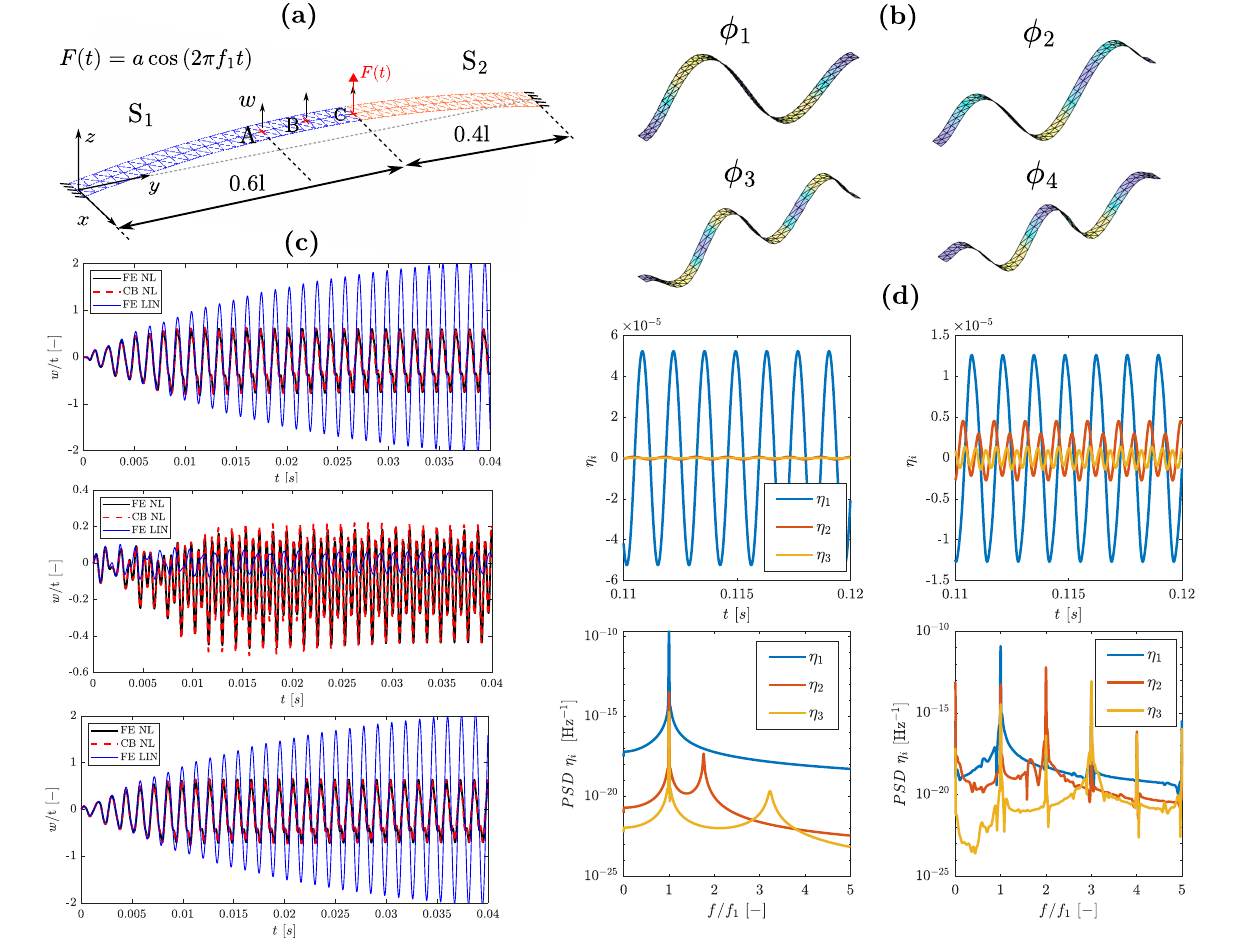}
    \caption{Response of a curved, clamped-clamped beam, investigated using the NL-CB ROM. In (a), FE assembly of the beam, partitioned into substructures $S_1$ and $S_2$. In (b), first $4$ VMs of the structure. In (c), time history of displacements along the $z$ axis at points A, B, and C, to applied harmonic load at resonance at point C. The different time histories were obtained using the nonlinear FE model (FE NL), the NL-CB ROM (CB NL), and the linearized FE model (FE LIN). In (d), time history and spectrogram of the amplitudes of global VMs $\phi_1$,$\phi_2$, and $\phi_3$, computed using the linearized FE model (left), and the nonlinear FE model (right).}
    \label{fig:curvedBeamResp}
\end{figure}

As a second example we consider a shallow curved beam, derived by adding uniform curvature to the flat beam presented in Section \ref{flatbeam}.
Specifically, the midspan of the beam is lifted of $5\ \text{mm}$, corresponding to $5\ \%$ of the span length. It is well known that curved beams pose a more challenging model reduction problem because of the direct quadratic coupling between low frequency modes \cite{Shen2021,Vizzaccaro2021}.
The assembly of the curved beam and the first four VMs are shown in Fig. \ref{fig:curvedBeamResp}a  and  Fig. \ref{fig:curvedBeamResp}b, respectively. 
The first four natural frequencies are $f_1 = 728.1\ \text{Hz} $,\ $f_2 = 1288.2\ \text{Hz} $,\ $f_3 = 2372.8\ \text{Hz} $, and $f_4 = 2860.2\ \text{Hz} $.
Rayleigh damping coefficients $\alpha = 58.46$ and $\beta = 1.58 \times 10^{-6}$ were identified by imposing a modal damping ratio of $1 \%$ for the first two modes.
\newl
The substructuring NL-CB ROM was constructed using four FIMs for each of the two substructures.
Interface reduction, following the same approach described for the flat beam, was also employed.
This yielded a ROM with $11$ DoFs.
The percentage error on the natural frequencies of the first four
VMs was respectively $e_1 = 0.26 \%$, $e_2 = 0.20 \%$, $e_3 = 0.56 \%$, and $e_4 = 0.15 \%$.
To validate the NL-CB ROM, we investigated the response of the structure subjected to a nodal force applied at $60\%$ of the span, in the direction of the $z$ axis, varying as $F(t) = \cos{(2\pi f_1 t)}$.
The NL-CB ROM, the full order FE model and its linearized counterpart were time integrated using the Newmark-$\beta$ time integration scheme with a fixed time step of $2.29 \times 10^{-5}\ \text{s}$. 
In Fig. \ref{fig:curvedBeamResp}c, we compare the time history of displacements along the $z$ axes at three different locations on the beam:  point A  at $x = 0.4\text{l}$, point B at $x = 0.5\text{l}$, and point C at $x = 0.6\text{l}$.
The nonlinear ROM is in good agreement with the full FE model. 
From the comparison with the linear model, we noted that the nonlinear coupling between the first VMs results in a strong participation of the second VMs to the nonlinear response.
This is shown in Fig. \ref{fig:curvedBeamResp}d, where we plot the time history of the modal amplitudes of the first three VMs, along with their frequency content, for the linear and nonlinear responses.
In the nonlinear case, the second VM has a strong harmonic component oscillating at twice the excitation frequency, as shown in the nonlinear spectrogram in Fig. \ref{fig:curvedBeamResp}d bottom right.
\subsection{Nonlinear Flat Panel}
\label{sec:panel}
\begin{figure}[h]
    \centering
    \includegraphics[width=1\linewidth]{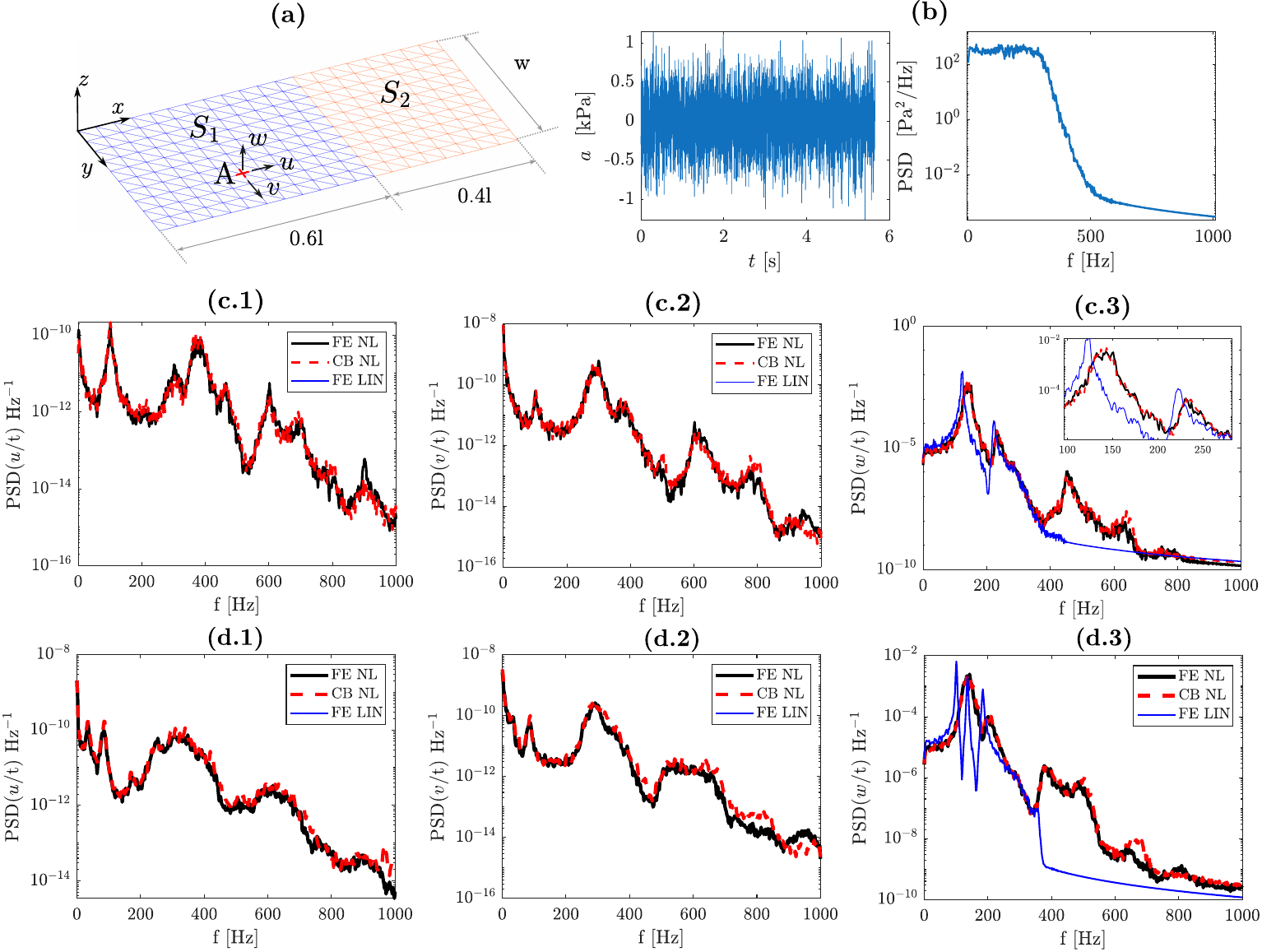}
    \caption{Nonlinear CB ROM used for investigating the response of a clamped flat panel to uniform acoustic pressure. In (a) FE assembly of the structure, partitioned into substructures $S_1$ and $S_2$. In (b) time variation of acoustic pressure (left) and its spectral content (right). In (c), PSDs of displacements time histories at point A along the $x$ (c.1), $y$ (c.2), and $z$(c.3) axes. The PSDs were computed using the nonlinear FE model (FE NL), the NL-CB ROM (CB NL), and a linearized FE model (FE LIN). In (d), same as in (c) but for modified thickness of substructure $S_2$. }
    \label{fig:flatPanelPSD}
\end{figure}

As third test case we consider a flat rectangular panel, clamped at all edges, of length $\text{l}=400\ \text{mm}$, width $\text{w} = 200 \ \text{mm}$, and thickness $\text{t} = 0.8\ \text{mm}$.
The material is linear elastic with an elastic modulus $E= 70\ \text{GPa}$, Poisson's ratio $\nu = 0.33$, and density $\rho = 2700 \ \text{Kg}/\text{m}^3$.
Geometric nonlinearities resulting from the von K\'arm\'an strain-displacement formulation were considered.
The panel was modeled using $400$ linear shell elements, each with $3$ nodes and $6$ DoFs per node, as shown in Fig. \ref{fig:flatPanelPSD}a.
The total number of DoFs, after applying the boundary conditions, amounted to $1026$.
The frequency of the first $10$ VMs is reported in table \ref{tab:panelFreqs}.
\begin{table}[h!]
\centering
\begin{tabular}{|c|c|c|c|c|c|c|c|c|c|c|}
\hline
mode & 1 & 2 & 3 & 4 & 5 & 6 & 7 & 8 & 9 & 10 \\ 
\hline
f [Hz] & 122.4 & 158.2 & 222.4 & 314.6  & 320.2  & 354.7  & 414.1 & 433.7 &  499.8 & 579.2  \\ 
\hline
$e_{\text{f}}$ $\%$ [-] & 0.05 & 0.02 & 0.10 & 0.28 & 0.28  & 0.15  & 0.24  & 0.22 & 0.82 &  3.39  \\ 
\hline
\end{tabular}
\caption{First $10$ natural frequencies of the clamped-clamped panel and their relative percentage error when computed using the substructuring ROM.}
\label{tab:panelFreqs}
\end{table}
Rayleigh damping with $\alpha = 8.672$ and $\beta = 1.134 \times 10^{-6}$ was considered. These values were obtained by imposing a damping ratio of $1\%$ to the first two VMs.
\newl
A partition of the structure in two substructures $S_1$ and $S_2$ was considered, as shown in Fig. \ref{fig:flatPanelPSD}a.
The interface was placed at $60\%$ of the span, parallel to the $y$ axes.
The NL-CB ROM was constructed using $8$ FIMs per substructure.
The highest frequency of the retained FIMs was $747.9\ \text{Hz}$ for $S_1$ and $1097.8\ \text{Hz}$ for $S_2$.
Interface reduction was used to reduce the number of interface Dofs, which amounts to 48.

The interface modes were computed by assemblying $S_1$ and $S_2$, computing the first $10$ VMs, and using their restriction to interface DoFs to form a interface RB with $10$ vectors.
Additionally, since the model was intended to be used in the nonlinear regime, complementary vectors were added to the interface RB to capture the nonlinear effects. 
Specifically, we appended to the interface RB the $6$ SMDs \cite{jain2017} defined from the first $3$ VMs and restricted to interface DoFs.
With this procedure we reduced the interface to $16$ DoFs.
The percentage error on the first $10$ natural frequencies obtained with the CB ROM is reported in Table \ref{tab:panelFreqs}.
\newl
To test the NL-CB ROM in the nonlinear regime, we considered the response of the panel subjected to acoustic pressure.
Specifically, we applied a uniform pressure in space, varying in time as a white noise with a sound pressure level of $144\ \text{dB}$ and a cut-off frequency of $300\ \text{Hz}$.
The time history of the pressure and its correspondent frequency spectrum are shown in Fig. \ref{fig:flatPanelPSD}b.
The NL-CB ROM, the full nonlinear FE model and its linearized counterpart were time integrated for a time span of $6\ \text{s}$ using the Newmark-$\beta$ time integration scheme with a fixed time step $h = 5.55\times 10^{-5}\ \text{s}$.
Power spectral densities (PSDs) of time histories of displacements for the probe node $\text{A}$, shown in Fig. \ref{fig:flatPanelPSD}a, were estimated using the Welch's method with time segments of $0.81\ \text{s}$ using $50\%$ overlap. 
As shown in Fig \ref{fig:flatPanelPSD}c, the PSDs obtained from displacements computed using the NL-CB ROM and the FE model are well overlaid, both for the in-plane (Fig. \ref{fig:flatPanelPSD}c.1 and Fig. \ref{fig:flatPanelPSD}c.2) and out-of-plane (Fig. \ref{fig:flatPanelPSD}c.3) motion.
Note that the linearized model  predicts zero in-plane displacements because in plane VMs are not directly excited by the load. 
Instead, they are only loaded through the bending-stretching coupling, as can be inferred from the nonlinear responses.
Moreover, from Fig. \ref{fig:flatPanelPSD}c.3, we observe that the nonlinearity produces lowering, smearing, and shifting to higher frequencies of the resonance peaks. 
Furthermore, it is worth noting the presence of an additional peak in the nonlinear response around $500\ \text{Hz}$, which is well beyond the load cut-off frequency.
\newl
The assembly procedure, required in the calculation of the VMs and SMDs for the computation of the interface RB, undermines the modularity of the substructuring approach. 
Other interface reduction techniques could be used.
Possible options are the System-Level Characteristic Constraint Modes reduction \cite{Castanier2001} or Local-Level Charachteristc Constraint Modes reduction with exact or weak interface compatibility \cite{Hong2013, Krattiger2019}.
The first approach is usually the most accurate, though it requires an assembly and the solution of an assembled eigenvalue problem. 
In contrast, in the second approach, the interface RB is assembled starting from the VMs of the single substructures.
However, also in these cases, if a substructure is modified, the interfaces between connected substructures must be recomputed and with it also the substructure ROMs.
To preserve modularity, one could think of using the same interface for different substructure configurations.
To test this concept, we modified substructure $S_2$ by reducing its thickness from $0.8\ \text{mm}$ to $0.5\ \text{mm}$ while keeping substructure $S_1$ unmodified.
Then, we updated the NL-CB ROM by recomputing only the ROM of substructure $S_2$, using the same interface RB that was computed for the uniform thickness panel.
Time integration was run for the same applied load.
The PSD of displacements at point A on the panel (see Fig. \ref{fig:flatPanelPSD}a) are plotted in Fig. \ref{fig:flatPanelPSD}d and show that minor accuracy is lost. This proves that the same interface can be potentially re-used to assess variations in the response to local substructure modifications. 
\subsection{MEMS Gyroscope}
\begin{figure}[H]
    \centering
    \includegraphics[width=1\linewidth]{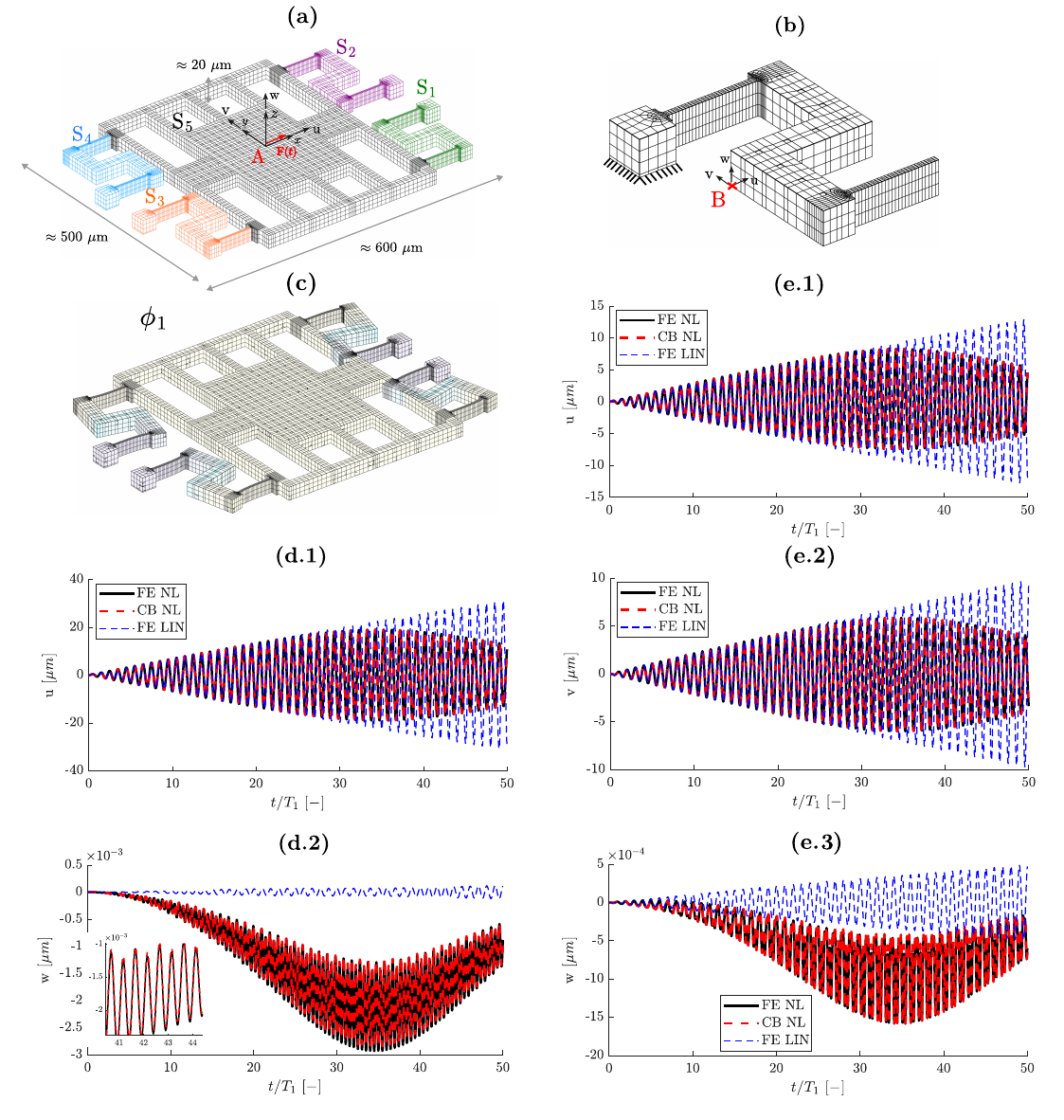}
    \caption{NL-CB ROM tested on a MEMS gyroscope. In (a) FE model of the MEMs gyroscope, partitioned into substructures $\text{S}_1,\text{S}_2,\text{S}_3,\text{S}_4,$ and $\text{S}_5$. In (b), detail of supporting beam. In (c), first VM. In (d) and (e) time history of displacements to harmonic excitation at resonance, computed using the nonlinear FE model (FE NL), the NL-CB ROM (CB NL) and the linearized FE model (FE LIN). In (d), displacements at point A along $x$ axis (d.1) and $z$ axis (d.2). In (e), displacements at point B, along $x$ axis (e.1), $y$ axis (e.2), and $z$ axis (e.3).}
    \label{fig:mems}
\end{figure}

As a last test case, we constructed the NL-CB ROM for the MEMS gyroscope presented in \cite{Marconi2021}, a device used to sense angular velocities. 
The MEMS meshed geometry is shown in Fig. \ref{fig:mems}a.
This device consists of an inertial mass connected to the ground with four s-shaped beams. 
During operational conditions, it is driven at the first natural frequency  by a pair of electrodes that exert a force along the $x$ axis, resulting in in-plane motion  (x-y plane in Fig. \ref{fig:mems}a). 
This motion generates an out-of-plane displacement along the z axis when an external angular velocity is aligned to the $y$ axis, as a result of the Coriolis force.
Through a direct measure of the out-of-plane displacement, the angular rate can be inferred.
\newl
As shown in Fig. \ref{fig:mems}, we considered a partition of the device into five different substructures, where
substructures $\text{S}_1,\text{S}_2,\text{S}_3,\text{and}\  \text{S}_4$ are the anchoring beams whereas substructure $\text{S}_5$ is the suspended mass. 
This domain partition can be used for efficient re-analysis: in case the suspended mass or the supporting beams are modified -- in an optimization loop, for instance -- only the corresponding substructures need to be updated while the rest remains untouched.
The geometry of the FE substructures was discretized with hexahedral elements with quadratic shape functions, each with $20$ DoFs. 
We adopted the \textit{Green-Lagrange} strain formulation and a linear elastic material with density $\rho = 2300 \ \text{Kg}/\text{m}^3$, Young's Modulus $E=148\ \text{Gpa}$, and Poisson's ratio $\nu=0.23$.
We used $1581$ elements for each of the four supporting beams and $8202$ for the central mass, yielding  a model with $14\,626$ elements.
The bottom face of each supporting beam was encastered by imposing zero displacements along the three axes, as schematically shown in Fig. \ref{fig:mems}.b.
After imposing the boundary conditions, the total number of free variables in the model was $261\,495$.
The natural frequency of the assembled FE model is $f_1 = 3.3418\times 10^4\  \text{Hz}$ and the corresponding VM is plotted in Fig. \ref{fig:mems}.c.
Note that this mode is in the $x-y$ plane.
Rayleigh damping with $\alpha=105$ and $\beta=0$ was considered. 
\newl
The NL-CB ROM was constructed using two FIMs for the center mass and one FIM for each of the supporting beams.
Moreover, the interface between the suspended mass and the beams was reduced using rigid body interface modes as in the \text{Virtual Node} method \cite{Krattiger2019}.
In this way, each interface was reduced to three translational modes along the $x$, $y$, and $z$ axes and three rigid body rotations around the center of each interface.
This resulted in $7$ and $26$ interface DoFs, for the supporting beams and for the center mass, respectively. 
With interface reduction, the number of DoFs of the assembled NL-CB ROM is $30$.
The NL-CB ROM was constructed using the \textit{Euler} cluster of \textit{ETH\ Z\"urich}, requesting $10$ cores and $7\ \text{Gb}$ of RAM per core.
With this computational power, the NL-CB ROM was constructed in $1\ \text{h}\ 58\ \text{min}$.
\newl
To assess accuracy of the NL-CB ROM, we considered the transient response of the structure subjected to harmonic excitation at resonance.
Specifically, two concentrated forces were applied to the top and bottom nodes in the center of the suspended mass (point A in Fig. \ref{fig:mems}.a.), each one varying in time as $\text{F}(t) = 25 \times 10^{-6} \sin{(2\pi f_1 t)} \ \text{N}$.
The FE model, a linearized FE model, and the NL-CB ROM were time integrated for $50$ cycles using the Newmark-$\beta$ time integration scheme with a constant time step $h=1/(40f_1) = 7.4811\times 10^{-7} \ \text{s}.$
The computed time history of nodal displacements at the center of the suspended mass (point A in Fig. \ref{fig:mems}.a) along the $x$ and $y$ axes are plotted in Fig. \ref{fig:mems}.d.1 and Fig. \ref{fig:mems}.d.2, respectively.
Additionally, we plot the  displacement time histories of point B (see Fig. \ref{fig:mems}.b) along the $x$, $y$ and $z$ axes in Fig. \ref{fig:mems}.e.1, Fig. \ref{fig:mems}.e.2, and Fig. \ref{fig:mems}.e.3, respectively.
The displacements obtained using the NL-CB ROM overlaps the ones obtained using nonlinear FE model.
Instead, the linearized model grows indefinitely, as expected for  linear systems excited at resonance. 
As for the nonlinear response, it is interesting to note the modulation of the oscillations probably due to the presence of a frequency component close to the excitation frequency.
This modulation does not appear in the linear system because only one frequency component exists, as the transient occurs at the first natural frequency, which equals the excitation frequency.
Interestingly, although four order of magnitudes smaller than the in-plane response, the nonlinear out-of-plane response exhibits a strong mean component. 
It is difficult to assess the influence of the out-of-plane motion on the overall in-plane motion of the nonlinear structure. 
However, also these small displacements are well captured by the NL-CB ROM, as shown in Fig. \ref{fig:mems}.d.2 and Fig. \ref{fig:mems}.e.3.
\newl
We report the computational times for time-integrating the nonlinear FE model and the NL-CB ROM.
The FE model was time integrated using the \text{Euler} cluster of \textit{ETH Z\"urich}, using five cores 
with $10\ \text{GB}$ RAM per core. 
The total wall-clock time was $60\ \text{h}\ 46\ \text{min}$.
Instead, the ROM was time-integrated on a local machine equipped with $32\ \text{GB}$ RAM and processor Intel Core Ultra 7 255H, operating at $ 5.1\ \text{GHz}$.
The computational time for time integration of the ROM was $3.28 \ \text{s}$, corresponding to a speed up of $6.65\times 10^4$.
Furthermore, the speed-up considering also the computational overhead for constructing the NL-CB ROM was $30.89$, which still represents a remarkable computational advantage over the direct use of the FE model.
\section{Conclusion} \label{sec:Conclusion}
In this paper, we presented a substructuring ROM for geometrically nonlinear structures that leverages the polynomial structure of elastic forces arising from the nonlinear strain-displacement formulation based on the von Karman and Green-Lagrange strain tensor.
Reduction  is performed at substructure level and the assembled ROM is obtained using a primal assembly procedure.
Substructure reduction is achieved by approximating internal DoFs using a quadratic displacement manifold.
Specifically, the dynamic equilibrium of internal DoFs is written in a basis of FIMs, divided into low and high frequency modes.
Then, the inertia and damping of high frequency FIMs are neglected, resulting in  nonlinear static equations that are solved at second order using Taylor series. 
In this way, high frequency FIMs are statically enslaved to low frequency FIMs and interface DoFs. 
The latter can also be potentially reduced with interface reduction techniques.
The substructure ROM equations are obtained by forcing to zero the  projection of the FE substructure equations onto the manifold tangent space.
In the projection step, three approximations are made: (i) the seventh order polynomial reduced forces, obtained by projecting the FE elastic forces, are truncated at third order; (ii) the convective inertial and damping terms, arising from the quadratic manifold terms, are neglected; (iii) the projection of the external load onto the quadratic part of the manifold is assumed to be negligible.
This results in a structure-preserving ROM that is highly efficient for time integration.
\newl
The substructure manifold is computed by repeatedly solving linear systems of equations featuring the inversion of the same linear operator. 
Furthermore, the nonlinear reduced tensors are computed using the element level-projection procedure whereby the element nonlinear tensors are projected on the manifold operators restricted to element variables.
In this way, the assembly of full order tensors is avoided, allowing for feasible ROM construction also in the case of large dimensional FE models.
\newl
Notably, if the substructure NL-CB ROM is truncated at linear order, the classic CB substructuring ROM \cite{bampton1968} is retrieved. Thus the presented technique can be regarded as an extension of the classic CB method to geometrically nonlinear structures.
\newl
The numerical tests presented in the paper show the excellent approximation quality of the NL-CB ROM, both for small dimensional academic benchmarks and for high-dimensional industrial problems. 
The ROM delivers large speed-ups in time integration compared to the nonlinear FE model.
Furthermore, ROM construction is feasible with moderate computational resources.
\newl
To conclude, two improvement areas can be identified.
Specifically, the ROM construction requires partial intrusiveness in the FE source code for the computation of the reduced tensors, limiting its use to specialized FE programs. A possible remedy can be found in the indirect identification schemes for nonlinear reduced tensors \cite{Mignolet2013}, which would require some modifications towards their integration with the presented method.
Finally, a second critical aspect is represented by the need for efficient interface reduction. As shown in the panel test case in Section \ref{sec:panel}, interface reduction might require a full assembly of the structures, thus partially compromising modularity. 
Nevertheless, the NL-CB ROM proves to be an accurate and efficient substructuring ROM that could be integrated with fully modular interface reduction techniques that are yet to be devised.
\section*{Funding Sources}

The second author gratefully acknowledges the Air Force Office for Scientific Research (AFOSR) for the supporting fund "Fully Coupled Reduced Order Models for Aero-Thermo-Elastic Analysis of Hypersonic Airframes: FullCoRe", award number FA8655-22-1-7040.
\appendix
\section{Second Order Manifold Derivation}
\label{app:secOrdS}
The coefficients of the Taylor expansion of the manifold in Eq. \eqref{eq:manFast}, up to second order, are obtained by solving Eq. \eqref{eq:staticEqModes}, for $\vvm{\hat{\eta}}{}{}$, using perturbation theory.
After inserting the polynomial expression of the nonlinear elastic forces (Eq. \eqref{eq:NlTensSubs}) into Eq. \eqref{eq:staticEqModes}, we get
\begin{align}
\label{eq:expForc}
    \vvm{\hat{\Phi}}{}{T} & \vv{K}{uu}{}  \vvm{\hat{\Phi}}{}{} \vvm{\hat{\eta}}{}{} + \vvm{\hat{\Phi}}{}{T} \vv{K}{uq}{} \vvm{\Psi}{}{}\vvm{\chi}{}{} + \vvm{\hat{\Phi}}{}{T} \vv{K}{uuu}{2} : (\vvm{\Phi}{}{} \vvm{\eta}{}{} \otimes \vvm{\Phi}{}{}\vvm{\eta}{}{}) + \vvm{\hat{\Phi}}{}{T} \vv{K}{uuu}{2} :(\vvm{\Phi}{}{}\vvm{\eta}{}{} \otimes \vvm{\hat{\Phi}}{}{}\vvm{\hat{\eta}}{}{}) +\nonumber \\
    & + \vvm{\hat{\Phi}}{}{T} \vv{K}{uuu}{2} :(\vvm{\hat{\Phi}}{}{}\vvm{\hat{\eta}}{}{} \otimes \vvm{\Phi}{}{}\vvm{\eta}{}{}) + \vvm{\hat{\Phi}}{}{T} \vv{K}{uuu}{2} :(\vvm{\hat{\Phi}}{}{}\vvm{\hat{\eta}}{}{} \otimes \vvm{\hat{\Phi}}{}{}\vvm{\hat{\eta}}{}{}) + \vvm{\hat{\Phi}}{}{T} \vv{K}{uuq}{2} :(\vvm{\Phi}{}{}\vvm{\eta}{}{} \otimes \vvm{\Psi}{}{}\vvm{\chi}{}{})+\nonumber \\
    &+ \vvm{\hat{\Phi}}{}{T} \vv{K}{uuq}{2} :(\vvm{\hat{\Phi}}{}{}\vvm{\hat{\eta}}{}{} \otimes \vvm{\Psi}{}{}\vvm{\chi}{}{}) + \vvm{\hat{\Phi}}{}{T} \vv{K}{uqq}{2} :(\vvm{\Psi}{}{}\vvm{\chi}{}{} \otimes \vvm{\Psi}{}{}\vvm{\chi}{}{}) +\mathcal{O}(3) = \zero.
\end{align}
 Then, we consider the Taylor expansion of $\vvm{\hat{\eta}}{}{}$ and explicitly write all the terms up to second order, obtaining 
 \begin{equation}
     \vvm{\hat{\eta}}{}{} = \sum_{i=1}^{\ap{n}{\phi}{}} \vv{\hat{A}}{i}{}\eta_i + \sum_{i=1}^{\ap{n}{\chi}{}} \vv{\hat{B}}{i}{}\chi_i+ \sum_{i=1}^{\ap{n}{\phi}{}}\sum_{j=i}^{\ap{n}{\phi}{}} \vv{\hat{C}}{ij}{} \eta_i \eta_j + \sum_{i=1}^{\ap{n}{\phi}{}}\sum_{j=1}^{\ap{n}{\chi}{}} \vv{\hat{D}}{ij}{} \eta_i \chi_j + \sum_{i=1}^{\ap{n}{\chi}{}}\sum_{j=i}^{\ap{n}{\chi}{}} \vv{\hat{E}}{ij}{} \chi_i \chi_j
     + \mathcal{O}(3),
     \label{eq:ansatz}
 \end{equation}
 where $\eta_i$ and $\chi_i$ are generic entries in vectors $\vvm{\eta}{}{}$ and $\vvm{\chi}{}{}$, whereas 
 \begin{equation}
     \vv{\hat{A}}{i}{} = \vv{\hat{A}}{}{}(:,i), \ \vv{\hat{B}}{i}{} = \vv{\hat{B}}{}{}(:,i),\  \vv{\hat{C}}{ij}{} = \vv{\hat{C}}{}{}(:,i,j),\  \vv{\hat{D}}{ij}{} = \vv{\hat{D}}{}{}(:,i,j) \in \R{\ap{n}{\hat{\phi}}{}}
 \end{equation} are the vector components of the multidimensional matrices defined in Eq. \eqref{eq:manFast}.
 For instance, using index notation, $[\hat{\text{C}}_{ij}]_k = \hat{\text{C}}_{kij} $.
 By inserting Eq. \eqref{eq:ansatz} in Eq. \eqref{eq:expForc},
 we get a polynomial expansion of the residual force in $\eta_i$ and $\chi_i$.
 Next, we collect all the terms multiplying the same monomials in $\eta_i$ and $\chi_i$, up to second order, and equate them to zero. By doing so, we obtain  a set of linear equations for the unknown ansatz coefficients.
At leading order we get
\begin{subequations}
 \label{eq:pertFirstOrd}
    \begin{align}
       \mathcal{O}(\eta_i) : \quad \vvm{\hat{\Phi}}{}{T} \vv{K}{uu}{} \vvm{\hat{\Phi}}{}{} \vv{\hat{A}}{i}{} & =  \zero\\
         \mathcal{O}(\chi_i) : \quad \vvm{\hat{\Phi}}{}{T} \vv{K}{uu}{} \vvm{\hat{\Phi}}{}{} \vv{\hat{B}}{i}{} & =  - \vvm{\hat{\Phi}}{}{T} \vv{K}{uq}{} \vvm{\Psi}{i}{},
         \label{eq:coeffb}
    \end{align}
\end{subequations}
while at second order
\begin{subequations}
\label{eq:pertSecondOrd}
    \begin{align}
    \mathcal{O}(\eta_i\eta_j) : \nonumber \\
       \vvm{\hat{\Phi}}{}{T} \vv{K}{uu}{} \vvm{\hat{\Phi}}{}{} \vv{\hat{C}}{ij}{} = &  \frac{1}{2}
          (2-\delta_{ij}) (- \vvm{\hat{\Phi}}{}{T} \vv{K}{uuu}{2} : ( \vvm{\Phi}{i}{} \otimes \vvm{\Phi}{j}{}) -\vvm{\hat{\Phi}}{}{T} \vv{K}{uuu}{2} : (\vvm{\Phi}{j}{} \otimes \vvm{\Phi}{i}{})),\label{eq:mdpr1}\\
        \mathcal{O}(\eta_i\chi_j) : \nonumber \\ \vvm{\hat{\Phi}}{}{T} \vv{K}{uu}{} \vvm{\hat{\Phi}}{}{} \vv{\hat{D}}{ij}{} = & - \vvm{\hat{\Phi}}{}{T} \vv{K}{uuu}{2}: (\vvm{\Phi}{i}{} \otimes \vvm{\hat{\Phi}}{}{}\vv{\hat{B}}{j}{}) -  \vvm{\hat{\Phi}}{}{T} \vv{K}{uuu}{2}: (\vvm{\hat{\Phi}}{}{}\vv{\hat{B}}{j}{} \otimes \vvm{\Phi}{i}{}) - \vvm{\hat{\Phi}}{}{T} \vv{K}{uuq}{2}: (\vvm{\Phi}{i}{} \otimes \vvm{\Psi}{j}{}),\\
        \mathcal{O}(\chi_i\chi_j) :\ \nonumber \\ 
          \vvm{\hat{\Phi}}{}{T} \vv{K}{uu}{} \vvm{\hat{\Phi}}{}{} \vv{\hat{E}}{ij}{} = & 
          \frac{1}{2}
          (2-\delta_{ij})(- \vvm{\hat{\Phi}}{}{T} \vv{K}{uuu}{2}: (\vvm{\hat{\Phi}}{}{}\vv{\hat{B}}{i}{} \otimes \vvm{\hat{\Phi}}{}{}\vv{\hat{B}}{j}{})  - \vvm{\hat{\Phi}}{}{T} \vv{K}{uuu}{2}: (\vvm{\hat{\Phi}}{}{}\vv{\hat{B}}{j}{} \otimes \vvm{\hat{\Phi}}{}{}\vv{\hat{B}}{i}{}) + \nonumber \\  & -\vvm{\hat{\Phi}}{}{T} \vv{K}{uuq}{2}: (\vvm{\hat{\Phi}}{}{} \vv{\hat{B}}{i}{}\otimes\vvm{\Psi}{j}{})  -
          \vvm{\hat{\Phi}}{}{T} \vv{K}{uuq}{2} : (\vvm{\hat{\Phi}}{}{} \vv{\hat{B}}{j}{}\otimes \vvm{\Psi}{i}{}) + \nonumber \\ & -\vvm{\hat{\Phi}}{}{T} \vv{K}{uqq}{2} :(\vvm{\Psi}{i}{} \otimes \vvm{\Psi}{j}{}) - \vvm{\hat{\Phi}}{}{T} \vv{K}{uqq}{2} :(\vvm{\Psi}{j}{} \otimes \vvm{\Psi}{i}{}))
    \end{align}
\end{subequations}
Here, $\delta_{ij}$ is the Kronecker delta and we denote with $\vvm{\Phi}{i}{},\vvm{\Psi}{i}{}$ the $i$th column vector in $\vvm{\Phi}{}{}$ and $\vvm{\Psi}{}{}$, respectively.
This procedure enforces Eq. \eqref{eq:staticEqModes} to be satisfied up to second order.
As expected in perturbation problems, the second order coefficients, $\vv{\hat{C}}{ij}{},\vv{\hat{D}}{ij}{},\  \text{and}\ \vv{\hat{E}}{ij}{}$, depend on the first order ones, $\vv{\hat{A}}{i}{}\ \text{and}\  \vv{\hat{B}}{i}{}$. Therefore, Eqs. \eqref{eq:pertSecondOrd} must be solved after Eqs. \eqref{eq:pertFirstOrd}. 
Additionally, the solution to these equations requires inversion of the same linear operator $ \vvm{\hat{\Phi}}{}{T} \vv{K}{uu}{} \vvm{\hat{\Phi}}{}{}$.
However, as detailed in the following, the manifold coefficients are not computed by directly using Eqs. \eqref{eq:pertFirstOrd} and \eqref{eq:pertSecondOrd}, as this would be computationally intractable for large FE systems. Instead, the alternative computational procedure presented in Section \ref{sec:manComp} is followed.
\section{Relation with Craig Bampton Reduction}
\label{sec:equivCB}
In this Section, we prove that if we truncate the reduction manifold expansion at linear level we retrieve the reduction space used in classical CB reduction.
\newl
The linear manifold operator, defined in Section \ref{sec:manDer}, writes
\begin{equation}
\label{eq:linearOpMan1}
    \vv{L}{\Gamma}{} = 
    \begin{bmatrix}
        \vvm{\Phi}{}{} & \vvm{\hat{\Phi}}{}{} \vv{\hat{B}}{}{} \\ \zero & \vvm{\Psi}{}{}
    \end{bmatrix},
\end{equation}
where $ \vvm{\Phi}{}{} $ are low frequency FIMs, $\vvm{\hat{\Phi}}{}{}$ are high frequency FIMs, $\vvm{\Psi}{}{}$ is the interface RB, and matrix $\vv{\hat{B}}{}{}$ satisfies 
\begin{equation}
 \quad \vvm{\hat{\Phi}}{}{T} \vv{K}{uu}{} \vvm{\hat{\Phi}}{}{} \vv{\hat{B}}{}{}  =  - \vvm{\hat{\Phi}}{}{T} \vv{K}{uq}{} \vvm{\Psi}{}{},      \label{eq:coeffb1}
\end{equation}
as shown in \ref{app:secOrdS}.
This RB is equivalent to the RB used in CB (see Eq. \eqref{eq:CBbasisIntRed}), where, instead of $\vvm{\hat{\Phi}}{}{} \vv{\hat{B}}{}{}$, the SMs, $\vv{S}{}{} $, are used. \newl
The definition of SMs is
\begin{equation}
     \quad  \vv{S}{}{}  =  -  \vv{K}{uu}{-1} \vv{K}{uq}{} \vvm{\Psi}{}{}.     \label{eq:coeffb11}
\end{equation}
It can be shown that
\begin{equation}
    \vvm{\hat{\Phi}}{}{} \vv{\hat{B}}{}{} = (\vv{I}{}{}-\vvm{\Phi}{}{}\vvm{\Phi}{}{T}\vv{M}{uu}{})\vv{S}{}{}.
    \label{eq:smodes1}
\end{equation}
In fact, by substituting Eq. \eqref{eq:smodes1} into Eq. \eqref{eq:coeffb1}, using Eq. \eqref{eq:coeffb11} and stiffness orthogonality of FIMs, an identity is retrieved.
As a result, the reduction space spanned by $[\vvm{\Phi}{}{}, \vvm{\hat{\Phi}}{}{} \vv{\hat{B}}{}{} ]$ is equivalent to the reduction space spanned by $[\vvm{\Phi}{}{}, \vv{S}{}{} ]$, as $\vv{S}{}{}-\vvm{\hat{\Phi}}{}{} \vv{\hat{B}}{}{}$   lies in the subspace spanned by $\vvm{\Phi}{}{}$.
This proves the equivalence between the linear subspace spanned by $\vv{L}{\Gamma}{}$ and the CB reduction subspace.
\section{Derivation of ROM Nonlinear Operators
}
\label{sec:redNlTens}
In this Section, we derive the ROM operators associated with the reduced nonlinear elastic forces.
We truncate the reduced forces to third order, based on the ROM simplifying hypothesis detailed in Section \ref{sec:projSubsEq}.
Once more, in this Section we omit denoting with superscript $"(s)"$ substructure related quantities.
\newl
We start the derivation by considering the projection of the elastic forces on the manifold.
 The projection of the elastic forces writes
\begin{equation}
    \pder{\vvm{\Gamma}{d}{T}}{\vvm{\xi}{}{}}(\vv{K}{}{}\vvm{\Gamma}{d}{} + \vv{f}{}{}(\vvm{\Gamma}{d}{})) = 
     \pder{\vvm{\Gamma}{d}{T}}{\vvm{\xi}{}{}}\left(
     \vv{K}{}{}\vvm{\Gamma}{d}{} + \vv{K}{}{2}:(\vvm{\Gamma}{d}{} \otimes\vvm{\Gamma}{d}{}) + \vv{K}{}{3}:(\vvm{\Gamma}{d}{} \otimes \vvm{\Gamma}{d}{} \otimes \vvm{\Gamma}{d}{})\right),
     \label{eq:projForc1}
\end{equation}
where we used the tensorial representation of forces in Eq. \eqref{eq:elastForcPolSubs}. Note that the dependence of the manifold on $\vvm{\xi}{}{}$ is implicitly assumed.
If we insert into Eq. \eqref{eq:projForc1} the manifold expression
\begin{equation}
     \vv{\Gamma}{d}{}(\vvm{\xi}{}{}) := \vv{L}{\Gamma}{}\vvm{\xi}{}{} + \vv{Q}{\Gamma}{}:(\vvm{\xi}{}{}\otimes \vvm{\xi}{}{}).
     \label{eq:manifGenerSubsB1}
\end{equation}
and the tangent projector expression
\begin{equation}
    \label{eq:projector1}
      \pder{\vvm{\Gamma}{d}{}}{\vvm{\xi}{}{}}(\vvm{\xi}{}{}) = \vv{L}{\Gamma}{} + \vv{P}{\Gamma}{}(\vvm{\xi}{}{}) =\vv{L}{\Gamma}{} + \vv{Q}{\Gamma,\tau}{} \cdot \vvm{\xi}{}{}.
\end{equation}
we obtain a polynomial equation of seventh order in $\vvm{\xi}{}{}$.
To derive the ROM, we compute the terms in the expansions only up to third order. 
\newl
For the following derivation, it is convenient to use the partition of the vector of reduced coordinates defined in Eq. \eqref{eq:decomposXi}:
\begin{equation}
    \vvm{\xi}{}{} = [\vvm{\eta}{}{T},\vvm{\chi}{}{T}]^T.
\end{equation}
Thus, the tangent projector derivative is rewritten using Eq. \eqref{eq:manifGenerSubsA} as
\begin{equation}
\label{eq:projExp}
     \pder{\vvm{\Gamma}{d}{T}}{\vvm{\xi}{}{}} =
     \begin{bmatrix}
         \pder{\vvm{\Gamma}{d}{T}}{\vvm{\eta}{}{}}\\
         \pder{\vvm{\Gamma}{d}{T}}{\vvm{\chi}{}{}}
     \end{bmatrix}
 = 
     \begin{bmatrix}
         \vvm{\Phi}{}{T} + \pder{\vvm{\Gamma}{\etahat}{T}}{\vvm{\eta}{}{}} \vvm{\hat{\Phi}}{}{T} & \zero \\ \pder{\vvm{\Gamma}{\etahat}{T}}{\vvm{\chi}{}{}} \vvm{\hat{\Phi}}{}{T} & \vvm{\Psi}{}{T}
     \end{bmatrix}.
\end{equation}
The partition of the tangent projector into its constant and linear components, $\vv{L}{\Gamma}{}$ and $\vv{P}{\Gamma}{}(\vvm{\xi}{}{})$ can be obtained from Eq. \eqref{eq:projExp} and from the expression of the linear projector in Eq. \eqref{eq:linearOpMan}.
Thus, we get 
\begin{equation}
\label{eq:projExpPart}
 \pder{\vvm{\Gamma}{d}{T}}{\vvm{\xi}{}{}} =
 \underbrace{
     \begin{bmatrix}
         \vvm{\Phi}{}{T}  & \zero \\ \vv{B}{}{T} \vvm{\hat{\Phi}}{}{T} & \vvm{\Psi}{}{T}
    \end{bmatrix}}_{\vv{L}{\Gamma}{T}} +
    \underbrace{
     \begin{bmatrix}
         \pder{\vvm{\Gamma}{\etahat}{T}}{\vvm{\eta}{}{}} \vvm{\hat{\Phi}}{}{T} & \zero \\ ( \pder{\vvm{\Gamma}{\etahat}{T}}{\vvm{\chi}{}{}} - \vv{\hat{B}}{}{T} )\vvm{\hat{\Phi}}{}{T} & \zero
     \end{bmatrix}}_{\vv{P}{\Gamma}{T}}
\end{equation}
By inserting Eq. \eqref{eq:projector1} into Eq. \eqref{eq:manifGenerSubsB1}, the projected forces write
\begin{equation}
      \pder{\vvm{\Gamma}{d}{T}}{\vvm{\xi}{}{}}(\vv{K}{}{}\vvm{\Gamma}{d}{} + \vv{f}{}{}(\vvm{\Gamma}{d}{})) = \vv{L}{\Gamma}{T} (\vv{K}{}{}\vvm{\Gamma}{d}{} + \vv{f}{}{}(\vvm{\Gamma}{d}{})) + \vv{P}{\Gamma}{T} (\vv{K}{}{}\vvm{\Gamma}{d}{} + \vv{f}{}{}(\vvm{\Gamma}{d}{})).
\end{equation}
We now prove that the second term in this expression generates monomial terms in $\vvm{\xi}{}{}$ of order greater than three, or equivalently that 
\begin{equation}
      \vv{P}{\Gamma}{T} (\vv{K}{}{}\vvm{\Gamma}{d}{} + \vv{f}{}{}(\vvm{\Gamma}{d}{})) = \mathcal{O}(4).
\end{equation}
To prove this statement, we  start from substructure elastic forces  partitioned into internal and interface DoFs,
\begin{equation}
\vv{K}{}{}\vv{d}{}{} + \vv{f}{}{}(\vv{d}{}{}) = 
    \begin{bmatrix}
    \vv{K}{uu}{} & \vv{K}{uq}{} \\
   \vv{K}{qu}{} & \vv{K}{qq}{} 
\end{bmatrix}
\begin{bmatrix}
    \vv{u}{}{} \\ \vv{q}{}{} 
\end{bmatrix}
+
\begin{bmatrix}
    \vv{f}{u}{}(\vv{u}{}{},\vv{q}{}{}) \\ 
    \vv{f}{q}{}(\vv{u}{}{},\vv{q}{}{})
\end{bmatrix},
\end{equation}
and we evaluate these forces on the manifold,
\begin{equation}
\label{eq:manifGenerSubsA1}
    \vv{d}{}{} = 
    \begin{bmatrix}
        \vv{u}{}{} \\ \vv{q}{}{}
    \end{bmatrix} 
    \approx 
    \begin{bmatrix}
        \vvm{\Phi}{}{}\vvm{\eta}{}{} + \vvm{\hat{\Phi}}{}{} \vvm{\Gamma}{\etahat}{}(\vvm{\eta}{}{},\vvm{\chi}{}{}) \\ \vvm{\Psi}{}{} \vvm{\chi}{}{}
    \end{bmatrix}.
\end{equation}
 In this way, we obtain
\begin{equation}
\label{eq:forcManProjExp1}
    \vv{K}{}{}\vvm{\Gamma}{d}{} + \vv{f}{}{}(\vvm{\Gamma}{d}{}) = \begin{bmatrix}
        \vv{f}{u}{el}(\vvm{\eta}{}{},\vvm{\chi}{}{}) \\ \vv{f}{q}{el}(\vvm{\eta}{}{},\vvm{\chi}{}{})
    \end{bmatrix},
\end{equation}
with
\begin{subequations}
\label{eq:forcManProjExp2}
    \begin{align}
        \vv{f}{u}{el}(\vvm{\eta}{}{},\vvm{\chi}{}{}) &= \vv{K}{uu}{}( \vvm{\Phi}{}{}\vvm{\eta}{}{} + \vvm{\hat{\Phi}}{}{} \vvm{\Gamma}{\etahat}{}) + \vv{K}{uq}{} \vvm{\Psi}{}{} \vvm{\chi}{}{} +  \vv{f}{u}{}(\vvm{\hat{\Phi}}{}{}\vvm{\Gamma}{\etahat}{},\vvm{\Psi}{}{}\vvm{\chi}{}{}),\\
        \vv{f}{q}{el}(\vvm{\eta}{}{},\vvm{\chi}{}{}) &= \vv{K}{qu}{}( \vvm{\Phi}{}{}\vvm{\eta}{}{} + \vvm{\hat{\Phi}}{}{} \vvm{\Gamma}{\etahat}{}) + \vv{K}{qq}{} \vvm{\Psi}{}{} \vvm{\chi}{}{} +  \vv{f}{q}{}(\vvm{\hat{\Phi}}{}{}\vvm{\Gamma}{\etahat}{},\vvm{\Psi}{}{}\vvm{\chi}{}{}). 
    \end{align}
\end{subequations}
Using the expression for $\vv{P}{\Gamma}{}$ in Eq. \eqref{eq:projExpPart} and Eqs. \eqref{eq:forcManProjExp1} and \eqref{eq:forcManProjExp2}, we can write 
\begin{equation}
     \vv{P}{\Gamma}{T} (\vv{K}{}{}\vvm{\Gamma}{d}{} + \vv{f}{}{}(\vvm{\Gamma}{d}{})) = 
     \begin{bmatrix}
     \underbrace{
       \pder{\vvm{\Gamma}{\etahat}{T}}{\vvm{\eta}{}{}} }_{\mathcal{O}(1)} 
       \underbrace{\vvm{\hat{\Phi}}{}{T} \vv{f}{u}{el}(\vvm{\eta}{}{},\vvm{\chi}{}{})}_{\mathcal{O}(3)} \\
       \underbrace{(\pder{\vvm{\Gamma}{\etahat}{T}}{\vvm{\chi}{}{}} - \vv{\hat{B}}{}{T})}_{\mathcal{O}(1)}
       \underbrace{\vvm{\hat{\Phi}}{}{T} \vv{f}{u}{el} (\vvm{\eta}{}{},\vvm{\chi}{}{}) }_{\mathcal{O}(3)} 
     \end{bmatrix} = \mathcal{O}(4).
     \label{eq:proj1}
\end{equation}
Here, we highlight the minimum polynomial order of the different terms and show that the total contributions yield monomials at least of fourth order.
In fact, $\vv{P}{\Gamma}{T}$ is linear by definition in $\vvm{\eta}{}{}$ and $\vvm{\chi}{}{}$, and $\vvm{\hat{\Phi}}{}{T} \vv{f}{u}{el} (\vvm{\eta}{}{},\vvm{\chi}{}{})$ is zero up to second order, since the displacement manifold was obtained by nullifying all terms up to second order in Eq. \eqref{eq:staticEqModes}.
Note that Eq. \eqref{eq:staticEqModes} is equivalent to $\vv{f}{u}{el}(\vvm{\eta}{}{},\vvm{\chi}{}{}) = \zero$.
\newl
In this way, we proved that
\begin{equation}
      \pder{\vvm{\Gamma}{d}{T}}{\vvm{\xi}{}{}}(\vv{K}{}{}\vvm{\Gamma}{d}{} + \vv{f}{}{}(\vvm{\Gamma}{d}{})) = \vv{L}{\Gamma}{T} (\vv{K}{}{}\vvm{\Gamma}{d}{} + \vv{f}{}{}(\vvm{\Gamma}{d}{})) + \mathcal{O}(4).
\end{equation}
Therefore, using the polynomial form of the elastic forces, we can write 
\begin{equation}
    \pder{\vvm{\Gamma}{d}{T}}{\vvm{\xi}{}{}}(\vv{K}{}{}\vvm{\Gamma}{d}{} + \vv{f}{}{}(\vvm{\Gamma}{d}{})) = 
     \vv{L}{}{T}\left(
     \vv{K}{}{}\vvm{\Gamma}{d}{} + \vv{K}{}{2}:(\vvm{\Gamma}{d}{} \otimes\vvm{\Gamma}{d}{}) + \vv{K}{}{3}:(\vvm{\Gamma}{d}{} \otimes \vvm{\Gamma}{d}{} \otimes \vvm{\Gamma}{d}{})\right) +\mathcal{O}(4).
     \label{eq:exp11}
\end{equation}
By substitution of Eq. \eqref{eq:manifGenerSubsB1} in Eq. \eqref{eq:exp11} and by discarding the other monomials greater than third order, we get
\begin{equation}
      \pder{\vvm{\Gamma}{d}{T}}{\vvm{\xi}{}{}}(\vv{K}{}{}\vvm{\Gamma}{d}{} + \vv{f}{}{}(\vvm{\Gamma}{d}{})) = \vv{K}{r}{} \vvm{\xi}{}{} + \vv{K}{r}{2} : (\vvm{\xi}{}{} \otimes \vvm{\xi}{}{}) + \vv{K}{r}{3} : (\vvm{\xi}{}{} \otimes \vvm{\xi}{}{} \otimes \vvm{\xi}{}{})  +\mathcal{O}(4),
      \label{eq:projForcMan}
\end{equation}
where
\begin{subequations}
\label{eq:tensExpressROMapp}
    \begin{align}
         \vv{K}{r}{} &= \vv{L}{\Gamma}{T}  \vv{K}{}{} \vv{L}{\Gamma}{}, \\
         \vv{K}{r}{2} &=  \vv{L}{\Gamma}{T} \cdot ( ( \vv{K}{}{2} \dotpr{21} \lman) \dotpr{21} \lman),\\
         \vv{K}{r}{3} &= \vv{L}{\Gamma}{T} \cdot (( \vv{K}{}{2} \dotpr{21} \lman) \dotpr{21} \qman ) + \vv{L}{\Gamma}{T} \cdot (( \vv{K}{}{2} \dotpr{21} \qman) \dotpr{21} \lman ) + \vv{L}{\Gamma}{T} \cdot ((( \vv{K}{}{3} \dotpr{21} \lman) \dotpr{21} \lman ) \dotpr{21} \lman).
    \end{align}
\end{subequations}
Note that, in this equation, the second order term  $\vv{L}{\Gamma}{T} \cdot \vv{K}{}{} \cdot \qman $ is not present.
This is because this term is null, as proved in the following.
\newl
We start by noting that the second order tensor $\vv{Q}{\Gamma}{}$ arises from the nonlinear part of the displacement manifold $\vvm{\hat{\Phi}}{}{}\vvm{\Gamma}{\etahat}{}$ in Eq. \eqref{eq:manifGenerSubsA1}, and thus can be written in general form as $[\vvm{\hat{\Phi}}{}{T},\zero]^T \cdot \vv{Q}{\etahat}{}$, where $\vv{Q}{\etahat}{}$ is a second order tensor associated with $\vvm{\Gamma}{\etahat}{}$.
Using the expression of the linear projector $\vv{L}{\Gamma}{T}$ in Eq. \eqref{eq:projExpPart}, we can write 
\begin{equation}
    \vv{L}{}{T} \vv{K}{}{} \cdot \vv{Q}{}{} = \begin{bmatrix}
        \vvm{\Phi}{}{T} &  \zero \\ \vvm{\hat{\Phi}}{}{} \vv{B}{}{T} & \vvm{\Psi}{}{T}
    \end{bmatrix}
    \begin{bmatrix}
        \vv{K}{uu}{} & \vv{K}{uq}{} \\ \vv{K}{qu}{} & \vv{K}{qq}{}
    \end{bmatrix} 
    \begin{bmatrix}
        \vvm{\hat{\Phi}}{}{T} \\ \zero 
    \end{bmatrix}
    \cdot
    \vv{Q}{\etahat}{}  = 
    \begin{bmatrix}
        \zero \\ (\vv{B}{}{T} \vvm{\hat{\Phi}}{}{T} \vv{K}{uu}{} + \vvm{\Psi}{}{T} \vv{K}{qu}{})\vvm{\hat{\Phi}}{}{}
    \end{bmatrix} \cdot \vv{Q}{\etahat}{}.
\end{equation}
With the definition of $\vv{B}{}{}$ matrix (see Eq. \eqref{eq:pertFirstOrd}), it is easy to see that this term vanishes, proving the thesis.
\section{Tangent Stiffness Derivative}
With reference to Eqs. \eqref{eq:subsEq1} and \eqref{eq:NlTensSubs}, the vector of substructure elastic forces acting on internal DoFs is defined as
\begin{equation}
    \vv{f}{u}{el}(\vv{u}{}{},\vv{q}{}{}) = \vv{K}{uu}{}\vv{u}{}{} + \vv{K}{uq}{} \vv{q}{}{} + \vv{f}{u}{}(\vv{u}{}{},\vv{q}{}{}),
    \label{eq:forcesss}
\end{equation}
where $\vv{f}{u}{}(\vv{u}{}{},\vv{q}{}{})$ are the nonlinear forces.
By explicitly expanding the nonlinear forces to second order, Eq. \eqref{eq:forcesss} writes
\begin{equation}
    \vv{f}{u}{el}(\vv{u}{}{},\vv{q}{}{}) = \vv{K}{uu}{}\vv{u}{}{} + \vv{K}{uq}{} \vv{q}{}{} +  \vv{K}{uuu}{2}:\vv{u}{}{}\otimes \vv{u}{}{} + \vv{K}{uuq}{2}:\vv{u}{}{}\otimes \vv{q}{}{} + \vv{K}{uqq}{2}:\vv{q}{}{}\otimes \vv{q}{}{} + \mathcal{O}(3).
\end{equation}
The substructure tangent stiffness matrix, partitioned into internal and interface DoFs, is derived by differentiating Eq. \eqref{eq:forcesss} with respect to $\vv{u}{}{}$ and $\vv{q}{}{}$, obtaining
\label{sec:tangStiffRel}
\begin{subequations}
\label{eq:tangStiffPart}
    \begin{align}
        \pder{\vv{f}{u}{el}}{\vv{u}{}{}}(\vv{u}{}{},\vv{q}{}{}) &=  \vv{K}{uu}{} + \vv{K}{uuu}{2} \cdot_{21} \vv{u}{}{} + \vv{K}{uuu}{2} \cdot_{31} \vv{u}{}{}   + \vv{K}{uuq}{2}\cdot_{31}\vv{q}{}{} +  \mathcal{O}(|\vv{u}{}{}|^2,|\vv{u}{}{}||\vv{q}{}{}|,|\vv{q}{}{}|^2),\\
        \pder{\vv{f}{u}{el}}{\vv{q}{}{}}(\vv{u}{}{},\vv{q}{}{}) &= \vv{K}{uq}{} + \vv{K}{uqq}{2} \cdot_{21} \vv{q}{}{} + \vv{K}{uqq}{2} \cdot_{31} \vv{q}{}{}  + \vv{K}{uuq}{2}\cdot_{21}\vv{u}{}{} +  \mathcal{O}(|\vv{u}{}{}|^2,|\vv{u}{}{}||\vv{q}{}{}|,|\vv{q}{}{}|^2)
    \end{align}
\end{subequations}
From Eqs. \eqref{eq:tangStiffPart}, we derive the expression of the tangent stiffness derivatives matrices computed around the origin, which write
\begin{subequations}
    \begin{align}
        \pder{}{\epsilon}\pder{\vv{f}{u}{el}}{\vv{u}{}{}}(\epsilon \vv{v}{u}{},\zero) \bigg{|}_{\epsilon=0} &=  \vv{K}{uuu}{2} \cdot_{21} \vv{v}{u}{} + \vv{K}{uuu}{2} \cdot_{31} \vv{v}{u}{}, \\
         \pder{}{\epsilon} \pder{\vv{f}{u}{el}}{\vv{u}{}{}}(\zero,\epsilon\vv{v}{q}{})\bigg{|}_{\epsilon=0} &=  \vv{K}{uuq}{2}\cdot_{31}\vv{v}{q}{},\\
         \pder{}{\epsilon} \pder{\vv{f}{u}{el}}{\vv{q}{}{}}(\epsilon \vv{v}{u}{},\zero) \bigg{|}_{\epsilon=0} &= \vv{K}{uuq}{2}\cdot_{21}\vv{v}{u}{}, \\
         \pder{}{\epsilon} \pder{\vv{f}{u}{el}}{\vv{q}{}{}}(\zero,\epsilon \vv{v}{q}{})\bigg{|}_{\epsilon=0} &=  \vv{K}{uqq}{2} \cdot_{21} \vv{v}{q}{} + \vv{K}{uqq}{2} \cdot_{31} \vv{v}{q}{}.
    \end{align}
\end{subequations}
Through right multiplication of the previous formulas with generic vectors $\vv{v}{u}{},\vv{v}{q}{},\vv{w}{u}{},\ \text{and}\ \vv{w}{q}{}$ we get
\begin{subequations}
    \begin{align}
       \pder{}{\epsilon}\pder{\vv{f}{u}{el}}{\vv{u}{}{}}(\epsilon \vv{v}{u}{},\zero) \bigg{|}_{\epsilon=0}  \cdot \vv{w}{u}{} &= \left( \vv{K}{uuu}{2} \cdot_{21} \vv{v}{u}{} + \vv{K}{uuu}{2} \cdot_{31} \vv{v}{u}{}\right) \cdot \vv{w}{u}{}, \\
         \pder{}{\epsilon} \pder{\vv{f}{u}{el}}{\vv{u}{}{}}(\zero,\epsilon\vv{v}{q}{})\bigg{|}_{\epsilon=0} \cdot \vv{v}{u}{}&=  \left( \vv{K}{uuq}{2}\cdot_{31}\vv{v}{q}{}\right) \cdot \vv{v}{u}{},\\
         \pder{}{\epsilon} \pder{\vv{f}{u}{el}}{\vv{q}{}{}}(\epsilon \vv{v}{u}{},\zero) \bigg{|}_{\epsilon=0} \cdot \vv{v}{q}{} &= \left( \vv{K}{uuq}{2}\cdot_{21}\vv{v}{u}{} \right) \cdot \vv{v}{q}{}, \\
         \pder{}{\epsilon} \pder{\vv{f}{u}{el}}{\vv{q}{}{}}(\zero,\epsilon \vv{v}{q}{})\bigg{|}_{\epsilon=0} \cdot \vv{w}{q}{} &=  \left( \vv{K}{uqq}{2} \cdot_{21} \vv{v}{q}{} + \vv{K}{uqq}{2} \cdot_{31} \vv{v}{q}{} \right) \cdot \vv{w}{q}{}.
    \end{align}
\end{subequations}
These expressions are the formulas used in Eq. \eqref{eq:tangStiffDer} for the numerical computation of the manifold coefficients.
\section{Algorithm for Computation of $\vvm{f}{}{*}$}
\label{sec:appAlg}
\begin{algorithm}[H]
\caption{Right Hand Side for Second Order Coefficients}
\textbf{Input:} \texttt{SubsFeModel}\sscr{a}, $\vvm{\Phi}{}{}$,  $\vvm{\Phi}{}{}\vv{\hat{B}}{}{}, \vvm{\Psi}{}{}$ \sscr{b} \\ 
\textbf{Output}: $\vv{F}{}{*}$ \sscr{c} .
\label{alg:ROMconstrRHS}
\begin{algorithmic}[1] 
    \Statex \textit{(i) Assemble r.h.s. of Quadratic Terms in $\eta$}
    \State $\vv{F}{}{*} \gets [\ ]$
    \For{$i=1,...,n_{\phi}$}
    \State $\text{d}\mathbf{K} \gets \texttt{SubsFeModel.dKuu\_u}(\vvm{\Phi}{}{}(:,i)) $ \sscr{d}
    \For{$j=i,...,n_{\phi}$}
    \State $\vv{f}{}{*} \gets  -\text{d}\mathbf{K}\cdot \vvm{\Phi}{}{}(:,j)$
    \If{$i = j$}
    \State $\vv{f}{}{*} \gets \vv{f}{}{*}/2$ 
    \EndIf
    \State $\vv{F}{}{*} \gets [\vv{F}{}{*},\vv{f}{}{*}]$
    \EndFor
    \EndFor
    \Statex
     \Statex \textit{(ii) Assemble r.h.s. of Quadratic Terms in $\chi$}
     \For{$i=1,...,n_{\chi}$}
     \State $\text{d}\vv{K}{1}{} \gets \texttt{SubsFeModel.dKuu\_u}(\vvm{\hat{\Phi}}{}{}\vv{\hat{B}}{}{}(:,i))$ \sscr{d}
     \State $\text{d}\vv{K}{2}{} \gets \texttt{SubsFeModel.dKuq\_u}(\vvm{\hat{\Phi}}{}{}\vv{\hat{B}}{}{}(:,i))$ \sscr{e}
     \State $\text{d}\vv{K}{3}{} \gets \texttt{SubsFeModel.dKuu\_q}(\vvm{\Psi}{}{}(:,i))$ \sscr{f}
     \State $\text{d}\vv{K}{4}{} \gets \texttt{SubsFeModel.dKuq\_q}(\vvm{\Psi}{}{}(:,i))$ \sscr{g}
     \For{$j=i,...,n_{\chi}$}
     \State $\vv{f}{1}{*} \gets -\text{d}\vv{K}{1}{}\cdot \vvm{\hat{\Phi}}{}{}\vv{\hat{B}}{}{}(:,j)$, 
       $\vv{f}{2}{*} \gets - \text{d}\vv{K}{2}{}\cdot \vvm{\Psi}{}{}(:,j)$
       \State $\vv{f}{3}{*} \gets - \text{d}\vv{K}{3}{}\cdot \vvm{\hat{\Phi}}{}{}\vv{\hat{B}}{}{}(:,j)$, $\vv{f}{4}{*}  \gets - \text{d}\vv{K}{4}{}\cdot \vvm{\Psi}{}{}(:,j)$
       \State $\vv{f}{}{*} \gets \vv{f}{1}{*} + \vv{f}{2}{*}+ \vv{f}{3}{*}+ \vv{f}{4}{*}$
     \If{$i = j$}
    \State $\vv{f}{}{*} \gets \vv{f}{}{*}/2$ 
    \EndIf
     \State $\vv{F}{}{*} \gets [\vv{F}{}{*},\vv{f}{}{*}]$
     \EndFor
     \EndFor
     \Statex
    \Statex \textit{(ii) Assemble r.h.s. of Quadratic Terms in $\eta$ and $\chi$}
    \For{$i=1,...,n_{\eta}$}
    \State \text{d}$\vv{K}{1}{} \gets  \texttt{SubsFeModel.dKuu\_u}(\vvm{\Phi}{}{}(:,i))$
   \State \text{d}$\vv{K}{2}{} \gets  \texttt{SubsFeModel.dKuq\_u}(\vvm{\Phi}{}{}(:,i))$
    \For{$j=1,...,n_{\chi}$}
    \State $\vv{f}{}{*} \gets -\text{d}\vv{K}{1}{} \cdot \vvm{\hat{\Phi}}{}{} \vv{\hat{B}}{}{}(:,j) - \text{d}\vv{K}{2}{} \cdot \vvm{\Psi}{}{}(:,j) $
    \State $\vv{F}{}{*} \gets [\vv{F}{}{*},\vv{f}{}{*}]$
    \EndFor
    \EndFor
\end{algorithmic}
\end{algorithm}
\vspace{-10pt}
\footnotesize
\noindent 
\sscr{a}\texttt{SubsFeModel} is the substructure FE model.\\
\sscr{b} Interface RB.\\
\sscr{c} Matrix that stores the rhs forces needed in the computation of the quadratic  manifold coefficients, as defined by Eq.\eqref{eq:rhsSecOrd}.\\
\sscr{d} Function \texttt{SubsFeModel.dKuu\_u}($\vv{v}{u}{}$) returns the substructure tangent stiffness derivative restricted to partition $\vv{K}{uu}{}$, along the substructure displacement with $\vv{u}{}{} = \vv{v}{u}{}$ and $\vv{q}{}{} = \zero.$ \\
\sscr{e} Function \texttt{SubsFeModel.dKuq\_u}($\vv{v}{u}{}$) returns the substructure tangent stiffness derivative restricted to partition $\vv{K}{uq}{}$, along the substructure displacement with $\vv{u}{}{} = \vv{v}{u}{}$ and $\vv{q}{}{} = \zero.$ \\
\sscr{f} Function \texttt{SubsFeModel.dKuu\_q}($\vv{v}{q}{}$) returns the substructure tangent stiffness derivative, restricted to partition $\vv{K}{uu}{}$, along the substructure displacement with $\vv{u}{}{} = \zero$ and $\vv{q}{}{} = \vv{v}{q}{}.$ \\
\sscr{f} Function \texttt{SubsFeModel.dKuq\_q}($\vv{v}{q}{}$) returns the substructure tangent stiffness derivative, restricted to partition $\vv{K}{uq}{}$, along the substructure displacement with $\vv{u}{}{} = \zero$ and $\vv{q}{}{} = \vv{v}{q}{}.$ 
\normalsize
\section{Reduced Elastic Potential Energy}
\label{sec:redPotEn}
In this Appendix, we derive the elastic potential function underlying the substructure ROM presented in Section \ref{sec:method}.
\newl
Let us consider the substructure FE model potential function, $\mathcal{V}_{FE}(\vv{d}{}{}) \in \R{}$.
If we insert the displacement manifold approximation $\vv{d}{}{} \approx \vvm{\Gamma}{d}{}(\vvm{\xi}{}{}) $ in $\mathcal{V}_{FE}(\vv{d}{}{})$ and differentiate with respect to $\vvm{\xi}{}{}$, we get
\begin{equation}
    \pder{\ }{\vvm{\xi}{}{}}  \mathcal{V}_{FE}(\vvm{\Gamma}{\vv{d}{}{}}{}(\vvm{\xi}{}{}))^T  =  \pder{\vvm{\Gamma}{d}{T}}{\vvm{\xi}{}{}}\pder{ \mathcal{V}_{FE}^T}{\vv{d}{}{}}\bigg\rvert_{\vvm{\Gamma}{d}{}} =  \pder{\vvm{\Gamma}{d}{T}}{\vvm{\xi}{}{}} (\vv{K}{}{}\vvm{\Gamma}{d}{}(\vvm{\xi}{}{}) + \vv{f}{}{}(\vvm{\Gamma}{d}{}(\vvm{\xi}{}{}))),
\end{equation}
where we used that the derivative of the FE potential function with respect to the nodal variables yields the nodal forces.
Using Eq. \eqref{eq:projForcMan} we can write
\begin{equation}
      \pder{\ }{\vvm{\xi}{}{}}  \mathcal{V}_{FE}(\vvm{\Gamma}{\vv{d}{}{}}{}(\vvm{\xi}{}{}))^T = \pder{\vvm{\Gamma}{d}{T}}{\vvm{\xi}{}{}}(\vv{K}{}{}\vvm{\Gamma}{d}{} + \vv{f}{}{}(\vvm{\Gamma}{d}{})) = \vv{f}{r}{}(\vvm{\xi}{}{}) + \mathcal{O}(4),
      \label{eq:projForcMan1}
\end{equation}
where
\begin{equation}
    \vv{f}{r}{}(\vvm{\xi}{}{}) = \vv{K}{r}{} \vvm{\xi}{}{} + \vv{K}{r}{2} : (\vvm{\xi}{}{} \otimes \vvm{\xi}{}{}) + \vv{K}{r}{3} : (\vvm{\xi}{}{} \otimes \vvm{\xi}{}{} \otimes \vvm{\xi}{}{})
    \label{eq:forceredpot}
\end{equation}
are the substructure ROM elastic forces. 
The potential function $ \mathcal{V}_{FE}(\vvm{\Gamma}{\vv{d}{}{}}{}(\vvm{\xi}{}{}))$ is a polynomial of eighth order in $\vvm{\xi}{}{}$.
If $\mathcal{V}_{FE}(\vvm{\Gamma}{\vv{d}{}{}}{}(\vvm{\xi}{}{}))_4$ is the truncation of $\mathcal{V}_{FE}(\vvm{\Gamma}{\vv{d}{}{}}{}(\vvm{\xi}{}{}))$ to fourth order monomials, we can write
\begin{equation}
    \mathcal{V}_{FE}(\vvm{\Gamma}{\vv{d}{}{}}{}(\vvm{\xi}{}{})) = \mathcal{V}_{FE}(\vvm{\Gamma}{\vv{d}{}{}}{}(\vvm{\xi}{}{}))_4 +\mathcal{O}(5),
\end{equation}
and by differentiating with respect to $\vvm{\xi}{}{}$ we get
\begin{equation}
     \pder{\ }{\vvm{\xi}{}{}}  \mathcal{V}_{FE}(\vvm{\Gamma}{\vv{d}{}{}}{}(\vvm{\xi}{}{}))^T = \pder{\ }{\vvm{\xi}{}{}}  \mathcal{V}_{FE}(\vvm{\Gamma}{\vv{d}{}{}}{}(\vvm{\xi}{}{}))_4^T + \mathcal{O}(4).
\end{equation}
Then, by comparison of this equation with Eq. \eqref{eq:projForcMan1} we can infer that
\begin{equation}
    \vv{f}{r}{}(\vvm{\xi}{}{}) = \pder{\ }{\vvm{\xi}{}{}}  \mathcal{V}_{FE}(\vvm{\Gamma}{\vv{d}{}{}}{}(\vvm{\xi}{}{}))_4^T 
\end{equation}
and, as a result, $\mathcal{V}_{r}(\vvm{\xi}{}{}) := \mathcal{V}_{FE}(\vvm{\Gamma}{\vv{d}{}{}}{}(\vvm{\xi}{}{}))_4$ is the ROM elastic potential function.
Additionally,
\begin{equation}
    \mathcal{V}_r(\vvm{\xi}{}{}) =  \frac{1}{2}\vvm{\xi}{}{T}\vv{K}{r}{} \vvm{\xi}{}{} + \frac{1}{3} \vvm{\xi}{}{} \cdot \vv{K}{r}{2} : (\vvm{\xi}{}{} \otimes \vvm{\xi}{}{}) + \frac{1}{4} \vvm{\xi}{}{} \cdot \vv{K}{r}{3} : (\vvm{\xi}{}{} \otimes \vvm{\xi}{}{} \otimes \vvm{\xi}{}{}),
\end{equation}
as by differentiating the above expression with respect to $\vvm{\xi}{}{}$ one retrieves Eq. \eqref{eq:forceredpot}. 

\section*{Acknowledgments}
The authors thank Jacopo Marconi for providing the original finite element model of the MEMS model discussed in the application section. 

\bibliography{bibbliography}

\end{document}